\documentclass[10pt]{amsart}
\usepackage{amsmath,amsfonts,amscd,amssymb,graphicx,mathrsfs,eufrak}
\usepackage[dvips,all,arc,curve,color,frame]{xy}
\newxyColor{pink}{1.0 0.4 0.5}{rgb}{}
\usepackage[usenames]{color}

\newtheorem{defin}{Definition}[section]
\newtheorem{thm}[defin]{Theorem}

\newtheorem{lemma}[defin]{Lemma}

\newtheorem{example}[defin]{Example}
\newtheorem{remark}[defin]{Remark}

\newcommand{\Mod}[1]{\mathfrak{Mod}#1}

\newcommand{\stq}[1]{\underline{\Mod{#1}}}

\newcommand{\Sstq}[1]{\textsf{St}(\stq{R})}

\newcommand{\C}[1]{\mathbb{C}^{#1}}

\renewcommand{\c}[1]{\mathcal{#1}}

\newcommand{\A}{\mathscr{A}}

\newcommand{\cl}[2]{c_{#1 #2}}
\newcommand{\an}[2]{a_{#1 #2}}

\newcommand{\CLr}[2]{{\color{Red}{C_{#1 #2}}}}
\newcommand{\AN}[2]{A_{#1 #2}}

\newcommand{\ff}[2]{f_{#1, #2}}
\renewcommand{\gg}[2]{g_{#1, #2}}
\newcommand{\FF}[2]{F_{#1, #2}}
\newcommand{\GG}[2]{G_{#1, #2}}
\newcommand{\pr}[2]{{\color{Red}{p_{#1 #2}}}}

\newcommand{\w}[1]{\widetilde{#1}}
\newcommand{\e}{\varepsilon}
\renewcommand{\d}[1]{\delta_{#1}}

\renewcommand{\t}[1]{\textnormal{#1}}
\renewcommand{\tt}[1]{\mathtt{#1}}

\def\minus{\hbox{-}}

\begin{document}
\title{\textsc{Reconstruction Algebras of Type $D$ (II)}}
\author{Michael Wemyss}
\address{Michael Wemyss\\ Graduate School of Mathematics\\ Nagoya University\\ Chikusa-ku, Nagoya, 464-8602, Japan}
\email{wemyss.m@googlemail.com}
\begin{abstract}
This is the third in a series of papers which give an explicit description of the reconstruction algebra as a quiver with relations; these algebras arise naturally as geometric generalizations of preprojective algebras of extended Dynkin quivers.  This paper is the companion to \cite{Wemyss_reconstruct_D(i)} and deals with dihedral groups $G=\mathbb{D}_{n,q}$ which have rank 2 special CM modules.  We show that such reconstruction algebras are described by combining a preprojective algebra of type $\w{D}$ with some reconstruction algebra of type $A$. 
\end{abstract}
\maketitle
\parindent 20pt
\parskip 0pt

\tableofcontents

\section{Introduction}
Everybody loves the preprojective algebra of an extended Dynkin diagram.  Algebraically they provide us with a rich source of non-trivial but manageable examples on which to test the latest theory, whilst geometrically they encode lots of information about  the singularities $\C{2}/G$ when $G$ is a finite subgroup of $\t{SL}(2,\C{})$.   However for quotients by finite subgroups of $\t{GL}(2,\C{})$ which are not inside $\t{SL}(2,\C{})$, preprojective algebras are not quite so useful.

Fortunately in this case the role of the preprojective algebra is played by a new algebra (the reconstruction algebra) which is built in a similar but slightly modified way.  This new algebra is by definition the endomorphism ring of the special CM modules in the sense of Wunram \cite{Wunram_generalpaper}, and the recent main result of \cite{Wemyss_GL2} states that for any finite group $G\leq\t{GL}(2,\C{})$ the quiver of the reconstruction algebra determines and is determined by the dual graph of the minimal resolution of $\C{2}/G$, labelled with self-intersection numbers.  Note that when the group $G$ is inside $\t{SL}(2,\C{})$ the reconstruction algebra is precisely the preprojective algebra of the associated extended Dynkin diagram.  

However it should also be noted that in general the above statement is known only on the level of quivers; to give an algebra we need to add in the extra information of the relations and it is this problem which is addressed in this paper.  Although technical in nature, knowledge of the relations of the reconstruction algebra unveils relationships between quotients arising from vastly different group structures and consequently reveals explicit geometric structure that has until now remained unnoticed.  In this paper, which should be viewed as a companion to \cite{Wemyss_reconstruct_D(i)},  we deal with dihedral groups $\mathbb{D}_{n,q}$ inside $GL(2,\C{})$ for which $n<2q$. The companion deals with the case $n>2q$.  The main difference between the two is that here rank two special CM modules enter the picture whereas in the $n>2q$ case all special CM modules have rank one.

Using knowledge of the specials from both \cite{Iyama_Wemyss_specials} and \cite{Wemyss_reconstruct_D(i)} we are able to associate to the dual graph of the minimal commutative resolution of $\C{2}/\mathbb{D}_{n,q}$ a quiver with relations and then prove that it is isomorphic to the endomorphism ring of the special CM modules.   Like most things involving explicit presentations of algebras the proof is rather technical, however we are able to side-step many issues by using both knowledge of the AR theory and knowledge of the geometry; in particular the non-explicit methods from \cite{Wemyss_GL2} are used extensively.  

We note that the reconstruction algebras in this paper naturally split into two subfamilies but it should perhaps be emphasized that the two families are fairly similar, particularly when compared to \cite{Wemyss_reconstruct_D(i)}.  The $\nu=N-1$ family (for notation see Section 2) should be viewed as being very geometrically similar to the dihedral ordinary double points, whereas the examples in \emph{loc. cit.} should be viewed as being very toric.  The remaining $0<\nu<N-1$ family sits geometrically somewhere between these two extremes.

This paper is organized as follows - in Section 2 we define the groups $\mathbb{D}_{n,q}$ and recap some combinatorics and results from \cite{Wemyss_reconstruct_D(i)} and \cite{Iyama_Wemyss_specials} which are needed for this paper.  Section 3 deals with the case when less than the maximal number of rank 2 specials occurs, whereas Section 4 deals with the case involving the maximal number of rank 2 specials.  In Section 5 we illustrate the correspondence between the algebra and the geometry by comparing two examples, making precise some of the above remarks which are of a more philosophical nature.

Throughout when working with quivers we shall write $ab$ to mean \textbf{$a$ followed by $b$.}  We work over the ground field $\mathbb{C}$ but any algebraically closed field of characteristic zero will suffice.

The author was supported with a JSPS Postdoctoral Fellowship held at Nagoya University, and the author would like to thank the JSPS for funding this work and Nagoya University for kind hospitality.  Thanks also to Osamu Iyama and Alvaro Nolla de Celis for useful conversations.

\section{Dihedral groups and special CM modules}
In this paper, as in \cite{Wemyss_reconstruct_D(i)}, we follow the notation of Riemenschneider \cite{Riemenschneider_invarianten}.
\begin{defin}
For $1<q<n$ with $(n,q)=1$ define the group $\mathbb{D}_{n,q}$ to be
\[
\begin{array}{cc}
\mathbb{D}_{n,q}=\left\{ \begin{array}{cc} \langle \psi_{2q},
\tau, \varphi_{2(n-q)} \rangle& \mbox{if } n-q\equiv 1 \mbox{ mod
}2\\\langle \psi_{2q}, \tau\varphi_{4(n-q)} \rangle& \mbox{if }
n-q\equiv 0 \mbox{ mod }2
\end{array}\right.
\end{array}
\]
with the matrices
{\scriptsize{
\[
\begin{array}{ccc}
\psi_k=\left( \begin{array}{cc}\e_k & 0\\ 0& \e_k^{-1}
\end{array}\right) &\tau=\left( \begin{array}{cc}0 & \e_4\\ \e_4& 0
\end{array}\right)&\varphi_k=\left( \begin{array}{cc}\e_k & 0\\ 0& \e_k
\end{array}\right)
\end{array}
\]
}}where $\e_t$ is a primitive $t^{th}$ root of unity.
\end{defin}
The order of the group $\mathbb{D}_{n,q}$ is $4(n-q)q$. By \cite[2.11]{Brieskorn} the dual graph of the minimal resolution of $\C{2}/\mathbb{D}_{n,q}$ is
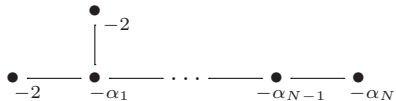
\[
\xymatrix@C=20pt@R=15pt{ &\bullet\ar@{-}[d]^<{-2}&&&\\
\bullet\ar@{-}[r]_<{-2} & \bullet\ar@{-}[r]_<{-\alpha_{1}}
&\hdots\ar@{-}[r] &\bullet\ar@{-}[r]_<{-\alpha_{N-1}} & \bullet
\ar@{}[r]_<{-\alpha_N}&}
\]
where the $\alpha$'s come from the Jung-Hirzebruch continued fraction expansion
\[
\frac{n}{q}=[\alpha_1,\hdots,\alpha_{N} ].
\]
By definition $\nu$ records the number of rank 2 special CM modules.  Denoting the dual continued fraction expansion by
\[
\frac{n}{n-q}=[a_2,\hdots,a_{e-1} ],
\]
we now briefly recap some combinatorics.  To the above data we define the series $c$, $d$, $r$, $i$, $l$, $b$, $\Delta$ and $\Gamma$ as follows:\\
1. The $c$ series, defined as $c_2=1$, $c_3=0$, $c_4=1$ and $c_{j}=a_{j-1}c_{j-1}-c_{j-2}$ for all $5\leq j \leq e$.\\
2.  The $d$ series, defined as $d_2=0$, $d_3=1$, $d_4=a_3-1$ and $d_{j}=a_{j-1}d_{j-1}-d_{j-2}$ for all $5\leq j \leq e$.\\
3.  The $r$ series, defined as $r_2=a_2(n-q)-q$, $r_3=r_2-(n-q)$, $r_4=(a_3+1)r_3-r_2$ and $r_j=a_{j-1}r_{j-1}-r_{j-2}$ for all $5\leq j\leq e$.\\
4. the $i$-series, defined as $i_0=n$, $i_1=q$ and $i_{t}=\alpha_{t-1}i_{t-1}-i_{t-2}$ for all $2\leq
t\leq N+1$.\\
5. The $l$-series, defined as $l_j=2+\sum_{p=1}^{j}(\alpha_p-2) \mbox{ for }1\leq j\leq N$.\\
6.  The $b$-series. Define $b_{0}:=1$, $b_{l_{N}-1}:=N$, and further for all  $1\leq t\leq l_N-2$ (if such $t$ exists), define $b_t$  to be the smallest integer $1\leq b_t\leq N$ such that $t\leq \sum_{p=1}^{b_t}(\alpha_p-2)$.\\
7. $\Delta_k$, defined as $\Delta_k=1+\sum_{t=\nu+1}^{k-1} c_{l_t} $ for all $\nu+1\leq k\leq N+1$.\\
8. $\Gamma_k$ defined as $\Gamma_k=\sum_{t=\nu+1}^{k-1} d_{l_t}$ for all $\nu+1\leq k\leq N+1$.\\ 
Note that in 7 and 8 we use the convention that for $k=\nu+1$ the sum is empty and so equals zero.

Now after setting $w_1=xy$ and
\[
\begin{array}{ll}
w_2=(x^q+y^q)(x^q+(-1)^{a_2}y^q)&w_3=(x^q-y^q)(x^q+(-1)^{a_2}y^q) \\
v_2=x^{2q}+(-1)^{a_2}y^{2q}&v_3=x^{2q}+(-1)^{a_2-1}y^{2q}
\end{array}
\]
we have the following result:
\begin{thm}\cite[Satz 2]{Riemenschneider_invarianten}\label{Riemen_generates} The
polynomials $w_1^{2(n-q)}$ and
$w_1^{r_t}v_2^{c_t}v_3^{d_t}$ for $2\leq t\leq e$ generate the
ring $\C{}[x,y]^{\mathbb{D}_{n,q}}$.  Alternatively we may take the polynomials $w_1^{2(n-q)}$ and
$w_1^{r_t}w_2^{c_t}w_3^{d_t}$ for $2\leq t\leq e$ as a generating set.
\end{thm}

The next two lemmas are the main combinatorial results which will be needed later to determine the relations on the reconstruction algebra:

\begin{lemma}\cite[2.12, 2.13]{Wemyss_reconstruct_D(i)}\label{r_is_difference_in_i_series}
Consider $\mathbb{D}_{n,q}$.  Then for all $2\leq t\leq e-2$, $r_{t+1}=r_{t+2}+i_{b_t}$.  Further $r_{l_t}=i_t-i_{t+1}$ for all $\nu+1\leq t\leq N$.
\end{lemma}
\begin{lemma}\cite[2.15]{Wemyss_reconstruct_D(i)}\label{c_and_d_lemma} Consider $\mathbb{D}_{n,q}$.  Then for all $2\leq t\leq e-2$,
\[
c_{t+2}=c_{t+1}+\Delta_{b_t}\quad\mbox{and}\quad d_{t+2}=d_{t+1}+\Gamma_{b_t}.
\]
In particular if $\nu=N-1$ then for all $3\leq t\leq e$
\[
c_{t}=t-3\quad\mbox{and}\quad d_{t}= 1.
\]
\end{lemma}

Define the rank 1 CM modules $W_+$, $W_-$ and for each $1\leq t\leq i_{\nu+1}+\nu (n-q)-1=q-1$ the rank 2 indecomposable CM module $V_t$ by the following positions in the
AR quiver of $\C{}[[x,y]]^{\mathbb{D}_{n,q}}$
{\tiny{
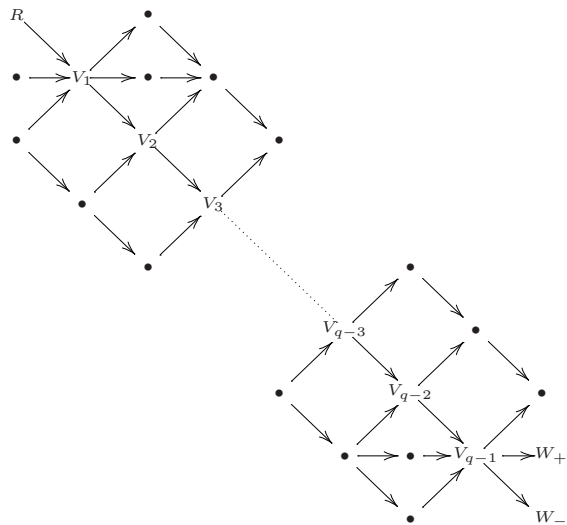
\[
\begin{array}{c}
\xymatrix@R=15pt@C=15pt{{}\save[]*{R}\restore\ar[dr(0.8)] &&\bullet\ar[dr(0.8)]&\\
\bullet\ar[r(0.8)]&{}\save[]*{V_1}\restore\ar[r(0.8)]\ar[ur(0.8)]\ar[dr(0.8)]&\bullet\ar[r(0.8)]&\bullet\ar[dr(0.8)]&\\
\bullet\ar[ur(0.8)]\ar[dr(0.8)]&&{}\save[]*{V_2}\restore\ar[ur(0.8)]\ar[dr(0.8)]&&\bullet&&\\
&\bullet\ar[ur(0.8)]\ar[dr(0.8)]&&{}\save[]*{V_3}\restore\ar[ur(0.8)]\ar@{.}[2,2]&&\\
&&\bullet\ar[ur(0.8)]&&&&\bullet\ar[dr(0.8)]&\\
&&&&&{}\save[]*{V_{q-3}}\restore\ar[dr(0.8)]\ar[ur(0.8)]&&\bullet\ar[dr(0.8)]&\\
&&&&\bullet\ar[dr(0.8)]\ar[ur(0.8)]&&{}\save[]*{V_{q-2}}\restore\ar[dr(0.8)]\ar[ur(0.8)]&&\bullet&\\
&&&&&\bullet\ar[r(0.8)]\ar[ur(0.8)]\ar[dr(0.8)]&\bullet\ar[r(0.6)]&{}\save[]*{V_{q-1}}\restore\ar[dr(0.8)]\ar[ur(0.8)]&{}\save[]*{\quad W_+}\restore\ar@{<-}[l(0.6)]\\
&&&&&&\bullet\ar[ur(0.8)]&&{}\save[]*{\quad W_-}\restore}
\end{array}
\]
}}i.e. all the $V_t$ lie on the diagonal leaving the vertex $R$, whilst $W_+$ and $W_-$ are the two rank 1 CM modules at the bottom of the diagonal.  Furthermore for every $1\leq t\leq n-q$ define the rank 1 CM module $W_t$ by the following position in the AR quiver:
{\tiny{
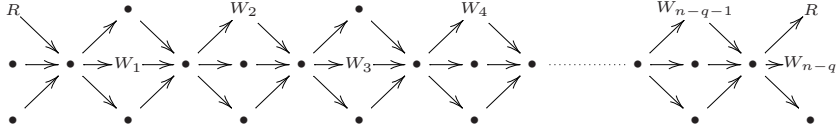
\[
\begin{array}{c}
\xymatrix@R=12pt@C=12pt{{}\save[]*{R}\restore\ar[1,1] &&\bullet\ar[1,1]&&{}\save[]*{W_2}\restore&&\bullet\ar[1,1]&&{}\save[]*{W_4}\restore&&&&{}\save[]*{W_{n-q-1}}\restore&&{}\save[]*{R}\restore\\
\bullet\ar[r]&\bullet\ar[r(0.7)]\ar[-1,1]\ar[1,1]&{}\save[]*{W_1}\restore&\bullet\ar@{<-}[l(0.75)]\ar[r]\ar[ur(0.8)]\ar[1,1]&\bullet\ar[r]&\bullet\ar[r(0.7)]\ar[-1,1]\ar[1,1]\ar@{<-}[ul(0.75)]&{}\save[]*{W_3}\restore&\bullet\ar@{<-}[l(0.75)]\ar[r]\ar[ur(0.8)]\ar[1,1]&\bullet\ar[r]&\bullet\ar@{<-}[ul(0.75)]\ar@{.}[0,2]&&\bullet\ar[r]\ar[ur(0.8)]\ar[1,1]&\bullet\ar[r]&\bullet\ar[r(0.5)]\ar[-1,1]\ar[1,1]\ar@{<-}[ul(0.75)]&{}\save[]*{W_{n-q}}\restore\\
\bullet\ar[-1,1]&&\bullet\ar[-1,1]&&\bullet\ar[-1,1]&&\bullet\ar[-1,1]&&\bullet\ar[-1,1]&&&&\bullet\ar[-1,1]&&\bullet}
\end{array}
\]
}}i.e. they all live on the non-zero zigzag leaving $R$.  Note that when $n-q$ is even the picture changes slightly since the position of $R$ on the right is twisted (for details we refer the reader to \cite[Section 6]{Iyama_Wemyss_specials}).  It is worth pointing out that $V_t=(\C{}[x,y]\otimes_{\C{}}\rho_t)^{\mathbb{D}_{n,q}}$ where $\mathbb{D}_{n,q}$ acts on both sides of the tensor (on the polynomial ring it acts as inverses) where $\rho_t$ is the two-dimensional irreducible representation%\footnote{just lift to $t^{th}$ powers since acts as inverses on the polynomial ring so cancels. Denote the basis of $V_t$ by $j_1^{(t)}$ and $j_2^{(t)}$ then $x^t\otimes j_1^{(t)}+y^t\otimes j_2^{(t)}\in V_t$}
\[
\begin{array}{c|c}
n-q \mbox{ odd} & n-q \mbox{ even}\\ \hline \begin{array}{ccl}
\psi_{2q}&\mapsto& \left(\begin{smallmatrix}\e_{2q}^{t}&0\\ 0&\e_{2q}^{-t}\end{smallmatrix}\right)
 \\ \tau & \mapsto &
\left(\begin{smallmatrix}0&\e_{4}^{t}\\ \e_{4}^{t}&0\end{smallmatrix}\right) \\ \varphi_{2(n-q)}&\mapsto &\left(\begin{smallmatrix}\e_{2(n-q)}^{t}&0\\ 0&\e_{2(n-q)}^{t}\end{smallmatrix}\right)
\end{array} & \begin{array}{ccl}
\psi_{2q}&\mapsto& \left(\begin{smallmatrix}\e_{2q}^{t}&0\\ 0&\e_{2q}^{-t}\end{smallmatrix}\right) \\\tau\varphi_{4(n-q)}&\mapsto &
\left(\begin{smallmatrix}0&\e_{4}^{t}\e_{4(n-q)}^{t}\\ \e_{4}^{t}\e_{4(n-q)}^{t}&0\end{smallmatrix}\right)
\end{array}
\end{array}
\]
For a description of the $W$'s in terms of representations, see \cite[Section 3]{Wemyss_reconstruct_D(i)}.

The next two results summarize the classification of the specials for the groups $\mathbb{D}_{n,q}$. 
\begin{thm}\cite[3.11]{Wemyss_reconstruct_D(i)}\label{1d_specials_generators}
For any $\mathbb{D}_{n,q}$, the following CM modules are special and
further they are 2-generated as
$\C{}[x,y]^{\mathbb{D}_{n,q}}$-modules by the following elements:
\[
\begin{array}{l|cc|cc}
W_{+} & x^q+y^q &(xy)^{n-q}(x^q-y^q)&x^q+y^q &(xy)^{n-q}(x^q-y^q) \\
W_{-} & x^q-y^q & (xy)^{n-q}(x^q+y^q)& x^q-y^q & (xy)^{n-q}(x^q+y^q) \\
W_{i_{\nu+1}} &(xy)^{i_{\nu+1}} & w_2=w_2^{\Delta_{\nu+1}}w_3^{\Gamma_{\nu+1}}&(xy)^{i_{\nu+1}} & v_2=v_2^{\Delta_{\nu+1}}v_3^{\Gamma_{\nu+1}}\\
W_{i_{\nu+2}} &(xy)^{i_{\nu+2}} & w_2^{\Delta_{\nu+2}}w_3^{\Gamma_{\nu+2}}&(xy)^{i_{\nu+2}} & v_2^{\Delta_{\nu+2}}v_3^{\Gamma_{\nu+2}}\\
&\vdots&&&\\
W_{i_N} & (xy)^{i_N} & w_2^{\Delta_{N}}w_3^{\Gamma_{N}}& (xy)^{i_N} & v_2^{\Delta_{N}}v_3^{\Gamma_{N}}
\end{array}
\]
where the left column is one such choice of generators, and the right hand column is another choice.
Further there are no other non-free rank one specials.
\end{thm}

\begin{thm}\cite[6.2]{Iyama_Wemyss_specials}
Consider the group $\mathbb{D}_{n,q}$ with $n<2q$,  then for all
$0\leq s\leq \nu-1$  $V_{i_{\nu+1}+s(n-q)}$ is special. Furthermore
these are all the rank 2 indecomposable special CM modules.
\end{thm}

For notational convenience denote $U_s:=V_{i_{\nu+1}+(\nu-s)(n-q)}$ for all $1\leq s\leq \nu$ , thus in this new notation the rank 2 indecomposable specials are simply $U_1,\hdots,U_\nu$.  Since for dihedral groups $\mathbb{D}_{n,q}$ the maximal co-efficient of any curve in the fundamental cycle $Z_f$ is two, by combining the previous two results we have full knowledge of all special CM modules.  We assign to them their corresponding vertices in the dual graph of the minimal resolution in Lemma~\ref{specials2vertices2} and Lemma~\ref{specials2vertices1} below.

Before proceeding to study the reconstruction algebra, we must first slightly re-interpret the preprojective algebra since this is crucial to later arguments.  Throughout the remainder of this section consider only the groups $G=\mathbb{D}_{n+1,n}\leq \t{SL}(2,\C{})$ and denote their natural representation by $V$ with basis $e_1$ and $e_2$.  It is very well known that the endomorphism ring of the specials in this case (summed without multiplicity) is isomorphic to the preprojective algebra of the extended Dynkin diagram of type $D$, and furthermore that this is morita equivalent to $\C{}[x,y]\# G$.  The relations for the preprojective algebra are obtained by choosing a $G$-equivariant basis for the maps in the McKay quiver; for details see \cite{CBH} or \cite{BSW}.  Since here for quivers we use the convention that $ab$ is $a$ followed by $b$, the number of arrows from $\rho$ to $\sigma$ in $\C{}[x,y]\# G$ is (see \cite[p6]{BSW} for details and notations)
\[
dim_{\C{}}e_{\rho}(V^*\otimes\C{}G)e_{\sigma}=dim_{\C{}}\t{Hom}_{\C{}G}(\C{}Ge_{\rho},V^*\otimes(\C{}Ge_{\sigma}))
\]
so we interpret arrows from vertex $\rho$ to vertex $\sigma$ in the preprojective algebra as $\C{}G$-maps from $\rho$ to $V^*\otimes\sigma$.  Since $G$ acts on $\C{}[x,y]$ in the same way as it acts on $V^*$, denote the basis of $V^*$ by $x$ and $y$.  The point is that the $\C{}G$-maps can be written in matrix form in terms of $x$'s and $y$'s, and composition of arrows is given by matrix multiplication.  To see this, for example denote the trivial representation of $\mathbb{D}_{n+1,n}$ by $\rho_0$ and denote its basis by $v_0$. Then we have a $\C{}G$-map
\[
\begin{array}{rcl}
\rho_0&\rightarrow &V^*\otimes V\\
v_0&\mapsto&x\otimes e_1+y\otimes e_2
\end{array}
\]
which we write as the matrix $\left(\begin{smallmatrix}x&y\end{smallmatrix}\right)$.  Similarly we have a $\C{}G$-map 
\[
\begin{array}{rcl}
V&\rightarrow&V^*\otimes \rho_0\\
e_1&\mapsto&(-y)\otimes v_0\\
e_2&\mapsto&x\otimes v_0
\end{array}
\]
which we write as the matrix $\left(\begin{smallmatrix}\minus y\\ x \end{smallmatrix}\right)$. The composition of these two maps is (see \cite[p7]{BSW})
\[
v_0\mapsto (x\otimes (-y)+y\otimes x)\otimes v_0
\]
which is zero in the skew group ring since $x\otimes y=y\otimes x$.  Thus we view the composition as being matrix multiplication $\left(\begin{smallmatrix}x&y\end{smallmatrix}\right)\left(\begin{smallmatrix}\minus y\\ x \end{smallmatrix}\right)$, subject to the relation that $x$ and $y$ commute.

We can choose the following $G$-basis for all the maps in the preprojective algebra
{\tiny
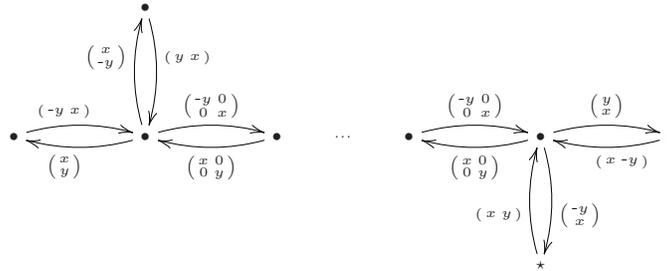
{
\[
\xymatrix@C=40pt@R=40pt{ &\bullet\ar@/^0.35pc/[d]^(0.4){\left(\begin{smallmatrix}
y&x \end{smallmatrix}\right)}&&&\\
\bullet\ar@/^0.35pc/[r]^(0.4){\left(\begin{smallmatrix}
\minus y&x \end{smallmatrix}\right)}&
\bullet\ar@/^0.35pc/[u]^(0.6){\left(\begin{smallmatrix}
x\\ \minus y \end{smallmatrix}\right)}\ar@/^0.35pc/[l]^(0.6){\left(\begin{smallmatrix}
x\\y \end{smallmatrix}\right)}\ar@/^0.35pc/[r]^{\left(\begin{smallmatrix}
\minus y&0\\0&x \end{smallmatrix}\right)}
&\bullet\ar@/^0.35pc/[l]^{\left(\begin{smallmatrix}
x&0\\0&y \end{smallmatrix}\right)}\ar@{}[r]|{\hdots} &\bullet\ar@/^0.35pc/[r]^{\left(\begin{smallmatrix}
\minus y&0\\0&x \end{smallmatrix}\right)}& \bullet\ar@/^0.35pc/[l]^{\left(\begin{smallmatrix}
x&0\\0&y \end{smallmatrix}\right)}\ar@/^0.35pc/[d]^(0.6){\left(\begin{smallmatrix}
\minus y\\ x \end{smallmatrix}\right)}\ar@/^0.35pc/[r]^{\left(\begin{smallmatrix}
y\\x \end{smallmatrix}\right)} & \bullet\ar@/^0.35pc/[l]^(0.4){\left(\begin{smallmatrix}
x&\minus y \end{smallmatrix}\right)}\\
&&&&\star\ar@/^0.35pc/[u]^(0.4){\left(\begin{smallmatrix}
x&y\end{smallmatrix}\right)}&}
\]
}}and thus we can view the preprojective algebra as being the above algebra subject to the rules of matrix multiplication and the relations induced by the fact that $x$ and $y$ commute. From this picture we may also read off the generators of the CM modules corresponding to the vertices.
\begin{remark}
\t{Note that the preprojective relations in the above picture at the two branching points are slightly different than the standard ones found elsewhere in the literature.}
\end{remark}
\begin{remark}
\t{The above picture is good for another reason, namely it gives us an explicit description of the Azumaya locus.  More precisely for any orbit containing the point $0\neq (\lambda_x,\lambda_y)\in\C{2}$ , denote by $M_{(\lambda_x,\lambda_y)}$ the representation of the preprojective algebra obtained by placing $\C{}$ at the outside vertices and $\C{2}$ at all other vertices, and further replacing every $x$ in the above picture by $\lambda_x$ and every $y$ by $\lambda_y$.  Then $M_{(\lambda_x,\lambda_y)}$ is the unique simple above the point corresponding to the orbit containing $(\lambda_x,\lambda_y)$.}
\end{remark}
%\footnote{For the change of basis isomorphism see handwritten sheet.  Note we are acting by inverses, and under this for example in $BD_{4.3}$ the change of basis under $\alpha$ is (reading left to right) $\left(\begin{smallmatrix}\e_8^3 &0\\0&\e_8^{-3}\end{smallmatrix}\right)$, $\left(\begin{smallmatrix}\e_8^2 &0\\0&\e_8^{-2}\end{smallmatrix}\right)$, $\left(\begin{smallmatrix}\e_8^1 &0\\0&\e_8^{-1}\end{smallmatrix}\right)$}

\section{The reconstruction algebra for $0<\nu<N-1$}
In this section we define the reconstruction algebra $D_{n,q}$ as a quiver with relations in the case $0<\nu<N-1$ and prove that it is isomorphic to the endomorphism ring of the special CM modules.   

Consider, for $N\in\mathbb{N}$ with $N\geq 2$ and for positive
integers $\alpha_{\nu+1}\geq 3$, $\alpha_t\geq 2$ for $\nu+2\leq t\leq N$, the labelled Dynkin diagram of type D:
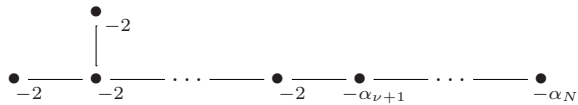
\[
\xymatrix@C=20pt@R=15pt{ &\bullet\ar@{-}[d]^<{-2}&&&\\
\bullet\ar@{-}[r]_<{-2} & \bullet\ar@{-}[r]_<{-2}
&\hdots\ar@{-}[r] &\bullet\ar@{-}[r]_<{-2} &\bullet\ar@{-}[r]_<{-\alpha_{\nu+1}}&\hdots\ar@{-}[r] & \bullet
\ar@{}[r]_<{-\alpha_N}&}
\]
Note that the hypothesis $0<\nu<N-1$ translates into the condition $\alpha_1=\hdots=\alpha_\nu=2$ with $\alpha_{\nu+1}\geq 3$. We call the left hand vertex the $+$ vertex, the top vertex the $-$ vertex and the remaining vertices $1,\hdots,N$ reading from left to right. 

To this picture we add an extended vertex called $\star$ and `double-up' as follows:
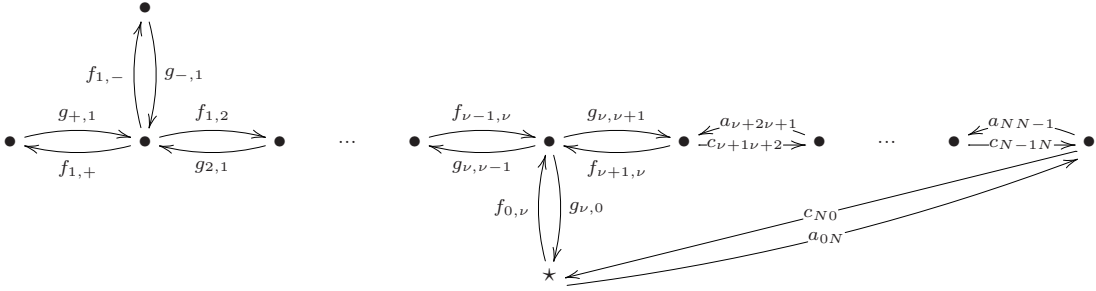
\[
\xymatrix@C=40pt@R=40pt{ &\bullet\ar@/^0.35pc/[d]^{\gg{-}{1}}&&&\\
\bullet\ar@/^0.35pc/[r]^{\gg{+}{1}}&
\bullet\ar@/^0.35pc/[u]^{\ff{1}{-}}\ar@/^0.35pc/[l]^{\ff{1}{+}}\ar@/^0.35pc/[r]^{\ff{1}{2}}
&\bullet\ar@/^0.35pc/[l]^{\gg{2}{1}}\ar@{}[r]|{\hdots} &\bullet\ar@/^0.35pc/[r]^{\ff{\nu-1}{\nu}}& \bullet\ar@/^0.35pc/[d]^{\gg{\nu}{0}}\ar@/^0.35pc/[l]^{\gg{\nu}{\nu-1}}\ar@/^0.35pc/[r]^{\gg{\nu}{\nu+1}} & \bullet\ar@/^0.35pc/[l]^{\ff{\nu+1}{\nu}}\ar@<-0.3ex>[r]|(0.45){\cl{\nu+1}{\nu+2}}&\bullet\ar@/_0.55pc/[l]|(0.45){\an{\nu+2}{\nu+1}}\ar@{}[r]|{\hdots} & \bullet\ar@<-0.3ex>[r]|{\cl{N-1}{N}} & \bullet\ar@/_0.55pc/[l]|{\an{N}{N-1}}\ar@<0.55ex>[1,-4]|{\cl{N}{0}} \\
&&&&\star\ar@/^0.35pc/[u]^{\ff{0}{\nu}}\ar@<-1.1ex>@/_0.55pc/[-1,4]|{\an{0}{N}}&&}
\]
where we have attached $\star$ to the $\nu^{th}$ and $N^{th}$ vertices.  Now if
$\sum_{i=1}^{N}(\alpha_i-2)\geq 2$ add extra arrows in the following way:
\begin{itemize}
\item if $\alpha_{\nu+1}>3$ then add an extra $\alpha_{\nu+1}-3$ arrows from the
$(\nu+1)^{st}$ vertex to $\star$.
\item If $\alpha_i>2$ with $i\geq \nu+2$ then add an extra $\alpha_i-2$ arrows from the $i^{th}$ vertex to $\star$.
\end{itemize}
Label the new arrows (if they exist) by $k_2,\hdots,k_{\sum
(\alpha_i-2)}$ reading left to right. Name this new quiver $Q$.%\footnote{Beginning with $k_2$ is not a typo!  Highest extra $k$ is $k_{e-3}$, after taking into account the $\cl{N}{0}$ the highest $k$ is $k_{e-2}$ (as in $\nu=0$)} 
\begin{example}\label{D138quiver}
\t{Consider $\mathbb{D}_{13,8}$ then $\frac{13}{8}=[2,3,3]$, $\nu=1$ and so the quiver $Q$ is}
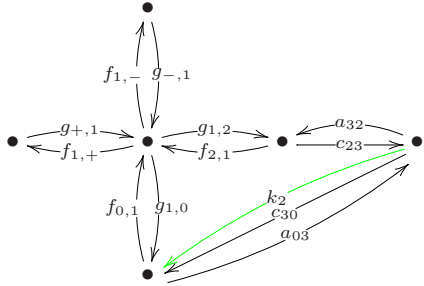
\[
\xymatrix@C=40pt@R=40pt{ &\bullet\ar@/^0.35pc/[d]|{\quad\gg{-}{1}}&&\\
\bullet\ar@/^0.35pc/[r]|{\gg{+}{1}}&
\bullet\ar@/^0.35pc/[u]|{\ff{1}{-}\quad}\ar@/^0.35pc/[l]|{\ff{1}{+}}\ar@/^0.35pc/[d]|{\quad\gg{1}{0}}\ar@/^0.35pc/[r]|{\gg{1}{2}}
&\bullet\ar@/^0.35pc/[l]|{\ff{2}{1}}\ar@<-0.3ex>[r]|{\cl{2}{3}} & \bullet\ar@/_0.55pc/[l]|{\an{3}{2}}\ar@[green]@<0.35ex>@/_0.5pc/[1,-2]|{k_2}\ar@<0.55ex>[1,-2]|{\cl{3}{0}}  \\
&\bullet\ar@/^0.35pc/[u]|{\ff{0}{1}\quad}\ar@<-1.1ex>@/_0.55pc/[-1,2]|{\an{0}{3}}&&}
\]
\end{example}
\begin{example}\label{D7356quiver}
\t{Consider $\mathbb{D}_{73,56}$ then $\frac{73}{56}=[2,2,2,5,2,3]$, $\nu=3$ and so the quiver $Q$ is}
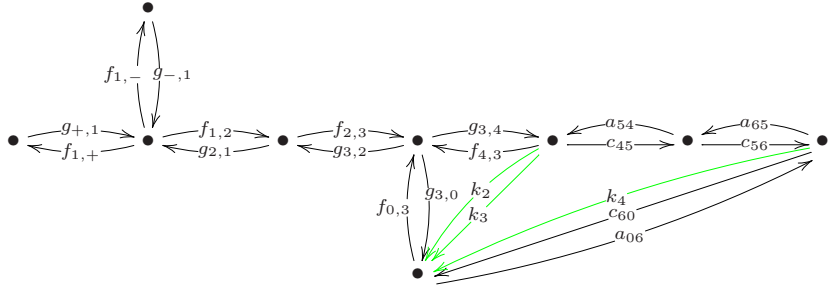
\[
\xymatrix@C=40pt@R=40pt{ &\bullet\ar@/^0.35pc/[d]|{\quad\gg{-}{1}}&&&\\
\bullet\ar@/^0.35pc/[r]|{\gg{+}{1}}&
\bullet\ar@/^0.35pc/[u]|{\ff{1}{-}\quad}\ar@/^0.35pc/[l]|{\ff{1}{+}}\ar@/^0.35pc/[r]|{\ff{1}{2}}&\bullet\ar@/^0.35pc/[r]|{\ff{2}{3}}\ar@/^0.35pc/[l]|{\gg{2}{1}}
&\bullet\ar@/^0.35pc/[d]|(0.4){\quad\gg{3}{0}}\ar@/^0.35pc/[l]|{\gg{3}{2}}\ar@/^0.35pc/[r]|{\gg{3}{4}} & \bullet\ar@/_0.5pc/@[green][1,-1]|(0.45){k_2}\ar@[green][1,-1]|(0.55){k_3}\ar@/^0.35pc/[l]|{\ff{4}{3}}\ar@<-0.3ex>[r]|{\cl{4}{5}} & \bullet\ar@/_0.55pc/[l]|{\an{5}{4}}\ar@<-0.3ex>[r]|{\cl{5}{6}} & \bullet\ar@/_0.55pc/[l]|{\an{6}{5}}\ar@<0.6ex>[1,-3]|{\cl{6}{0}}\ar@[green]@<0.4ex>@/_0.5pc/[1,-3]|{k_4} \\
&&&\bullet\ar@/^0.35pc/[u]|{\ff{0}{3}\quad}\ar@<-1.1ex>@/_0.55pc/[-1,3]|{\an{0}{6}}&&&}
\]
\end{example}

To define the reconstruction algebra we need to first set up some notations, then specify relations.   Firstly we denote $k_1:=\gg{\nu}{0}$ and $k_{1+\sum_{i=1}^{N}(\alpha_i-2)}:=\cl{N}{0}$.
 \begin{defin}
For all $1\leq r\leq 1+\sum_{i=1}^{N}(\alpha_i-2)$ define $\tt{B}_{r}$ (``the butt") to be the number of the vertex associated to the tail of the arrow $k_r$.
\end{defin}
Notice for all $2\leq r\leq 1+\sum_{i=1}^{N}(\alpha_i-2)$ it is true that $\tt{B}_r=b_r$ where $b_r$ is the $b$-series of $\frac{n}{q}$ defined in Section 2, however $\tt{B}_1\neq b_1$ since $b_1=\nu+1$ whilst $\tt{B}_1=\nu$ by the special definition of $k_1$.

Now define $u_{\nu}=1$ and further for $\nu+1\leq i\leq N$ denote
\[
\begin{array}{c}
u_i:=\t{max}\{ j: 2\leq j\leq 1+\sum_{i=1}^{N}(\alpha_i-2)\mbox{ with } b_{j}=i \}\\ 
v_i:=\t{min}\{ j: 2\leq j\leq 1+\sum_{i=1}^{N}(\alpha_i-2)\mbox{ with } b_{j}=i \}
\end{array}
\]
if such things exist (i.e. vertex $i$ has an extra arrow leaving it).  Also define $W_{\nu+1}:=\nu$ and for every $\nu+2\leq i\leq N$ define
\[
W_i=\left\{ \begin{array}{cl} \nu & \mbox{if }\c{S}_i\mbox{ is empty}\\ \mbox{the maximal number in }\c{S}_i & \mbox{else}
\end{array}\right. 
\]
where 
\[
\c{S}_i=\{\mbox{vertex }j: 1\leq j< i\mbox{ and }j\mbox{ has an extra arrow leaving it} \}.
\]
Thus $W_i$ is defined for $\nu+1\leq i\leq N$.  The idea behind it is that $W_i$ records the closest vertex to the left of vertex $i$ which has a $k$ leaving it; since we have defined $k_1:=\gg{\nu}{0}$ this is always possible to find.  Now define, for all $\nu+1\leq i\leq N$, $V_{i}:=u_{W_i}$.  Thus $V_i$ records the number of the largest $k$ to the left of the vertex $i$, where since $k_1:=\gg{\nu}{0}$ and $u_\nu=1$ it  always exists.

Denote also
\[
\FF{1}{\nu}:=\left\{  \begin{array}{cl} e_1&\mbox{if }\nu=1 \\ \ff{1}{2}\hdots\ff{\nu-1}{\nu}&\mbox{if }\nu>1 \end{array} \right.\qquad\GG{\nu}{1}:=\left\{  \begin{array}{cl} e_1&\mbox{if }\nu=1 \\ \gg{\nu}{\nu-1}\hdots\gg{2}{1}&\mbox{if }\nu>1 \end{array} \right.
\]
where $e_1$ is the trivial path at vertex $1$, and define $\AN{0}{t}:=\an{0}{N}\hdots\an{t+1}{t}$ for every $\nu+1\leq t\leq N$.

\begin{defin}\label{reconD2}
For $\frac{n}{q}=[\alpha_1,\hdots,\alpha_N]$ with $0<\nu<N-1$ define
$D_{n,q}$, the reconstruction algebra of type $D$, to be the path
algebra of the quiver $Q$ defined above subject to the relations
\[
\ff{0}{\nu}\gg{\nu}{\nu+1}=2\AN{0}{\nu}
\]
\[
\begin{array}{cl}
\mbox{If }\nu=1& \left\{\begin{array}{l}\gg{+}{1}\ff{1}{+}=0\quad \ff{0}{1}\gg{1}{0}=0 \\
\gg{-}{1}\ff{1}{-}=0 \quad \ff{2}{1}\gg{1}{2}=0 \\
\gg{1}{0}\ff{0}{1}-\gg{1}{2}\ff{2}{1}=\ff{1}{+}\gg{+}{1}-\ff{1}{-}\gg{-}{1}\end{array}\right.\\
&\\
\mbox{If }\nu>1& \left\{\begin{array}{l}\gg{+}{1}\ff{1}{+}=0\quad \ff{0}{\nu}\gg{\nu}{0}=0\\  
 \gg{-}{1}\ff{1}{-}=0\quad \ff{\nu+1}{\nu}\gg{\nu}{\nu+1}=0 \\
\ff{1}{+}\gg{+}{1}-\ff{1}{-}\gg{-}{1}=2\ff{1}{2}\gg{2}{1}\\
\gg{t}{t-1}\ff{t-1}{t}=\ff{t}{t+1}\gg{t+1}{t}\mbox{ for all }2\leq t\leq \nu-1\\
\gg{\nu}{0}\ff{0}{\nu}-\gg{\nu}{\nu+1}\ff{\nu+1}{\nu}=2\gg{\nu}{\nu-1}\ff{\nu-1}{\nu}
\end{array}\right.
\end{array}
\]
together with the relations defined algorithmically as 
\[
\begin{array}{ccl}
\mbox{Step }\nu+1&\mbox{If }\alpha_{\nu+1}=3 &\cl{\nu+1}{\nu+2}\an{\nu+2}{\nu+1}=\pr{\nu+1}{\nu+1}\\
&\mbox{If }\alpha_{\nu+1}>3 &k_2\AN{0}{\nu+1}=\pr{\nu+1}{\nu+1}, \AN{0}{\nu+1}k_2=\CLr{0}{\nu}k_1\\
&&k_t\CLr{0}{\nu+1}=k_{t+1}\AN{0}{\nu+1}, \CLr{0}{\nu+1}k_t=\AN{0}{\nu+1}k_{t+1}\,\,\forall\,\,2\leq t< u_{\nu+1}\\  & &
k_{u_{\nu+1}}\CLr{0}{\nu+1}=\cl{\nu+1}{\nu+2}\an{\nu+2}{\nu+1}\\  
&&\vdots\\
\mbox{Step }i&\mbox{If }\alpha_{i}=2  & \cl{i}{i+1}\an{i+1}{i}=\an{i}{i-1}\cl{i-1}{i}\\
 & \mbox{If }\alpha_i>2 & k_{v_i}\AN{0}{i}=\an{i}{i-1}\cl{i-1}{i}, \AN{0}{i}k_{v_i}=\CLr{0}{\tt{B}_{V_i}}k_{V_i}  \\
 & & k_{t}\CLr{0}{i}=k_{t+1}\AN{0}{i}, \CLr{0}{i}k_t=\AN{0}{i}k_{t+1}\,\,\forall\,\,v_i\leq t< u_i\\
 & & k_{u_i}\CLr{0}{i}=\cl{i}{i+1}\an{i+1}{i}\\ 
&&\vdots\\
\mbox{Step }N&\mbox{If }\alpha_{N}=2 & \cl{N}{0}\an{0}{N}=\an{N}{N-1}\cl{N-1}{N}, \an{0}{N}\cl{N}{0}=\CLr{0}{\tt{B}_{V_N}}k_{V_N}\\
&\mbox{If }\alpha_{N}>2 & k_{v_N}\an{0}{N}=\an{N}{N-1}\cl{N-1}{N},\an{0}{N}k_{v_N}=\CLr{0}{\tt{B}_{V_N}}k_{V_N}  \\
 & & k_{t}\CLr{0}{N}=k_{t+1}\an{0}{N}, \CLr{0}{N}k_{t}=\an{0}{N}k_{t+1}\,\,\forall\,\,v_N\leq t<
 u_N
\end{array}
\]
The only things that remain to be defined are the $\pr{}{}$ and the $\CLr{}{}$'s, and it is these which change according to the presentation, namely%\footnote{$\pr{\nu+1}{\nu+1}$ represents $(xy)^{r_3}w_2^{c_3}w_3^{d_3}=(xy)^{r_3}w_3$ in moduli presentation and $(xy)^{r_3}v_3$ in symmetric presentation.  They are the same iff $\nu$ is even.  The corresponding cycle at $\star$ is $\CLr{0}{\nu}k_1$, but note that in moduli presentation it goes to the opposite horn than $\pr{\nu+1}{\nu+1}$ does.  Note also that $\CLr{0}{\nu}k_1$ at vertex 0 is the same in both presentations iff $\nu$ is even.}
%\footnote{$\CLr{0}{\nu+1}$ represents $w_2$ in moduli presentation and $v_2$ in symmetric presentation.  They are the same iff $\nu$ is odd.}
%\footnote{Notice that its almost nice: i.e. that $\CLr{0}{\nu+1}=\CLr{0}{\nu}\gg{\nu}{\nu+1}$.  This is true in symmetric presentation, and also in moduli presentation when $\nu$ is odd, but is false in moduli presentation when $\nu$ is even.}
\begin{eqnarray*}
\pr{\nu+1}{\nu+1}&=&\left\{  \begin{array}{c} \ff{\nu+1}{\nu}\GG{\nu}{1}\ff{1}{+}\gg{+}{1}\FF{1}{\nu}\gg{\nu}{\nu+1}\\ \frac{1}{2}\ff{\nu+1}{\nu}\GG{\nu}{1}(\ff{1}{+}\gg{+}{1}+\ff{1}{-}\gg{-}{1})\FF{1}{\nu}\gg{\nu}{\nu+1} \end{array} \right.\\ 
\CLr{0}{\nu}&=&\left\{  \begin{array}{c} \ff{0}{\nu}\GG{\nu}{1}\ff{1}{-}\gg{-}{1}\FF{1}{\nu} \\ \frac{1}{2}\ff{0}{\nu}\GG{\nu}{1}(\ff{1}{+}\gg{+}{1}+\ff{1}{-}\gg{-}{1})\FF{1}{\nu} \end{array} \right.\\ 
\CLr{0}{\nu+1}&=&\left\{  \begin{array}{c}  \ff{0}{\nu}\GG{\nu}{1}\ff{1}{+}\gg{+}{1}\FF{1}{\nu}\gg{\nu}{\nu+1}\\ \frac{1}{2}\ff{0}{\nu}\GG{\nu}{1}(\ff{1}{+}\gg{+}{1}+\ff{1}{-}\gg{-}{1})\FF{1}{\nu}\gg{\nu}{\nu+1} \end{array} \right.
\end{eqnarray*}
where for the moduli presentation we take the top choices and for the symmetric presentation we take the bottom choices.  We also define, for all $\nu+2\leq t\leq N$, $\CLr{0}{t}=\CLr{0}{\nu+1}\cl{\nu+1}{\nu+2}\hdots\cl{t-1}{t}$.
\end{defin}
\begin{remark}
\t{As in \cite{Wemyss_reconstruct_D(i)} we should explain why we give two presentations.  The symmetric case is pleasing since it treats the two $(-2)$ horns equally, so that the algebra produced is independent of how we view the dual graph.  On the other hand the relations in the moduli presentation are more binomial and so this makes the explicit geometry easier to write down in Section 5.  We show in Theorem~\ref{iso_of_rings2} that the two presentations yield isomorphic algebras, but for the moment denote $D_{n,q}$ for the moduli presentation and $D_{n,q}^\prime$ for the symmetric presentation.}
\end{remark}
\begin{remark}
\t{The algorithmic relations are just the relations from the reconstruction algebra of Type $A$.  Thus the reconstruction algebra here should be viewed as a preprojective algebra of type $D$ with a reconstruction algebra of type $A$ stuck onto the side.}
\end{remark}
\begin{remark}
\t{The $\nu=1$ and $\nu>1$ split in the above definition is not natural in the sense that in both cases they are just the preprojective relations.  We only split them for notational ease.  In both presentations notice the existence of a relation which is not a cycle and further that at any vertex in the reconstruction algebra corresponding to a (-2) curve there is only one relation, a preprojective relation.}%\footnote{The symmetric and moduli presentations are \emph{exactly} the same if and only if we are in the base case and further $\nu$ is even (e.g. 0!).  To see this, note that $\nu$ even is required so that $\pr{\nu+1}{\nu+1}$ and $\CLr{0}{\nu}$ are the same in both presentations, whereas base case needed so that $\CLr{0}{\nu+1}$ doesn't enter the picture (since the two presentations of it are different for $\nu$ even)}
\end{remark}

\begin{example}\label{D138relations}
\t{For the group $\mathbb{D}_{13,8}$ the symmetric presentation of the reconstruction algebra is the quiver in Example~\ref{D138quiver} subject to the relations}
{\footnotesize{
\[
\begin{array}{c}
\ff{0}{1}\gg{1}{2}=2\an{0}{3}\an{3}{2}\\
\begin{array}{l}
\gg{+}{1}\ff{1}{+}=0 \\ \gg{-}{1}\ff{1}{-}=0\\ \ff{0}{1}\gg{1}{0}=0 \\
 \ff{2}{1}\gg{1}{2}=0\end{array}\\
\gg{1}{0}\ff{0}{1}-\gg{1}{2}\ff{2}{1}=\ff{1}{+}\gg{+}{1}-\ff{1}{-}\gg{-}{1}\\
\frac{1}{2}\ff{2}{1}(\ff{1}{+}\gg{+}{1}+\ff{1}{-}\gg{-}{1})\gg{1}{2}=\cl{2}{3}\an{3}{2}\\
\begin{array}{rl}
\an{3}{2}\cl{2}{3}=k_2(\an{0}{3})&
\frac{1}{2}\ff{0}{1}(\ff{1}{+}\gg{+}{1}+\ff{1}{-}\gg{-}{1})\gg{1}{0}=(\an{0}{3})k_2\\
k_2(\frac{1}{2}\ff{0}{1}(\ff{1}{+}\gg{+}{1}+\ff{1}{-}\gg{-}{1})\gg{1}{2}\cl{2}{3})=\cl{3}{0}(\an{0}{3})& (\frac{1}{2}\ff{0}{1}(\ff{1}{+}\gg{+}{1}+\ff{1}{-}\gg{-}{1})\gg{1}{2}\cl{2}{3})k_2=(\an{0}{3})\cl{3}{0}
\end{array}
\end{array}
\]}}
\end{example}
\begin{example}\label{D7356relations}
\t{For the group $\mathbb{D}_{73,56}$ the moduli presentation of the reconstruction algebra is the quiver in Example~\ref{D7356quiver} subject to the relations}
{\footnotesize{
\[
\begin{array}{c}
\ff{0}{3}\gg{3}{4}=2\an{0}{6}\an{6}{5}\an{5}{4}\\
\begin{array}{l}
\gg{+}{1}\ff{1}{+}=0\\\gg{-}{1}\ff{1}{-}=0 \\  
 \ff{0}{3}\gg{3}{0}=0\\ \ff{4}{3}\gg{4}{3}=0 
 \end{array} \\
 \begin{array}{c}
 \ff{1}{+}\gg{+}{1}-\ff{1}{-}\gg{-}{1}=2\ff{1}{2}\gg{2}{1}\\
\gg{2}{1}\ff{1}{2}=\ff{2}{3}\gg{3}{2}\\
\gg{3}{0}\ff{0}{3}-\gg{3}{4}\ff{4}{3}=2\gg{3}{2}\ff{2}{3}
\end{array}\\
\begin{array}{rl}
\ff{4}{3}\gg{3}{2}\gg{2}{1}\ff{1}{+}\gg{+}{1}\ff{1}{2}\ff{2}{3}\gg{3}{4}=k_2(\an{0}{6}\an{6}{5}\an{5}{4})&
\ff{0}{3}\gg{3}{2}\gg{2}{1}\ff{1}{-}\gg{-}{1}\ff{1}{2}\ff{2}{3}\gg{3}{0}=(\an{0}{6}\an{6}{5}\an{5}{4})k_2\\
k_2(\ff{0}{3}\gg{3}{2}\gg{2}{1}\ff{1}{+}\gg{+}{1}\ff{1}{2}\ff{2}{3}\gg{3}{4})=k_3(\an{0}{6}\an{6}{5}\an{5}{4})&
(\ff{0}{3}\gg{3}{2}\gg{2}{1}\ff{1}{+}\gg{+}{1}\ff{1}{2}\ff{2}{3}\gg{3}{4})k_2=(\an{0}{6}\an{6}{5}\an{5}{4})k_3
\end{array}\\
k_3(\ff{0}{3}\gg{3}{2}\gg{2}{1}\ff{1}{+}\gg{+}{1}\ff{1}{2}\ff{2}{3}\gg{3}{4})=\cl{4}{5}\an{5}{4}\\
\an{5}{4}\cl{4}{5}=\cl{5}{6}\an{6}{5}\\
\begin{array}{rl}
\an{6}{5}\cl{5}{6}=k_4(\an{0}{6})& (\an{0}{6})k_4=(\ff{0}{3}\gg{3}{2}\gg{2}{1}\ff{1}{+}\gg{+}{1}\ff{1}{2}\ff{2}{3}\gg{3}{4})k_3\\
k_4(\ff{0}{3}\gg{3}{2}\gg{2}{1}\ff{1}{+}\gg{+}{1}\ff{1}{2}\ff{2}{3}\gg{3}{4}\cl{4}{5}\cl{5}{6})=\cl{6}{0}(\an{0}{6})& (\an{0}{6})\cl{6}{0}=(\ff{0}{3}\gg{3}{2}\gg{2}{1}\ff{1}{+}\gg{+}{1}\ff{1}{2}\ff{2}{3}\gg{3}{4}\cl{4}{5}\cl{5}{6})k_4
\end{array}
\end{array}
\]}}
\end{example}
Before we show that the endomorphism ring of the specials is isomorphic to this quiver with relations we must first justify the positions of the specials relative to the dual graph:
\begin{lemma}\label{specials2vertices2}
Consider $\mathbb{D}_{n,q}$ with $0<\nu<N-1$.  Then the specials correspond to the dual graph of the minimal resolution in the following way
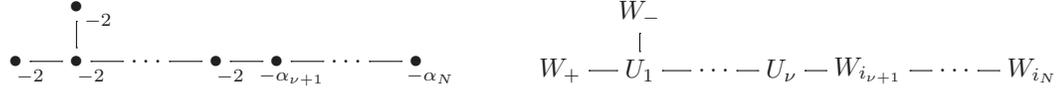
\[
\begin{array}{cc}
\begin{array}{c}
\xymatrix@C=12pt@R=10pt{ &\bullet\ar@{-}[d]^<{-2}&&&\\
\bullet\ar@{-}[r]_<{-2} & \bullet\ar@{-}[r]_<{-2}
&\hdots\ar@{-}[r] &\bullet\ar@{-}[r]_<{-2} &\bullet\ar@{-}[r]_<{-\alpha_{\nu+1}}&\hdots\ar@{-}[r] & \bullet
\ar@{}[r]_<{-\alpha_N}&}\end{array} &
\begin{array}{c}
\xymatrix@C=8pt@R=6pt{ &W_-\ar@{-}[d]&&&\\
W_+\ar@{-}[r] & U_1\ar@{-}[r]
&\hdots\ar@{-}[r] &U_\nu\ar@{-}[r] &W_{i_{\nu+1}}\ar@{-}[r]&\hdots\ar@{-}[r] & W_{i_N}}
\end{array}
\end{array}
\]
\end{lemma}
\begin{proof}
Denote by $E_t$ the curve corresponding to vertex $t$ for $t\in\{ +,-,1,\hdots,N \}$.  The assumption $0<\nu< N-1$ translates into the condition $\alpha_{\nu+1}\geq 3$ which then forces
$Z_f=E_{+}+E_{-}+\sum_{p=1}^{\nu}2E_p+\sum_{t=\nu+1}^{N}E_t$.  Hence
\[
\begin{array}{cc|l}
-Z_f\cdot E_t&(Z_K-Z_f)\cdot E_t&\\ \hline
0&0&t\in\{+,-,1,\hdots,\nu-1\}\\ 
1&1&t=\nu\\
\alpha_{\nu+1}-3&-1&t=\nu+1\\
\alpha_{t}-2&0&\nu+1<t<N\\
\alpha_N-1&1&t=N 
\end{array}
\]
and so by  \cite[3.3]{Wemyss_GL2} the quiver of the endomorphism ring of the specials coincides with $Q$ as defined above.  But it is clear that the following composition of maps in the AR quiver do not factor through other specials
{\tiny{
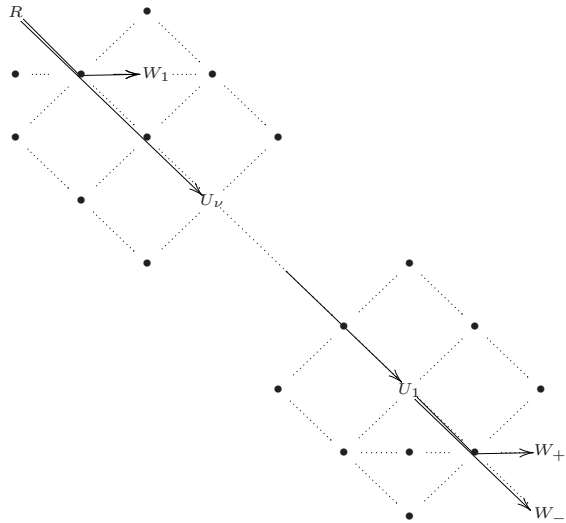
\[
\begin{array}{c}
\xymatrix@R=15pt@C=15pt{{}\save[]*{R}\restore\ar@<-0.5ex>[3,3]\ar`dr/0pt[1,1][1,2] &&\bullet\ar@{.}[dr(0.8)]&\\
\bullet\ar@{.}[r(0.5)]&\bullet\ar@{.}[ur(0.8)]\ar@{.}[dr(0.8)]&{}\save[]*{\quad W_1}\restore\ar@{<-}[l(0.5)]&\bullet\ar@{.}[l(0.6)]\ar@{.}[dr(0.8)]&\\
\bullet\ar@{.}[ur(0.8)]\ar@{.}[dr(0.8)]&&\bullet\ar@{.}[ur(0.8)]\ar@{.}[dr(0.8)]&&\bullet&&\\
&\bullet\ar@{.}[ur(0.8)]\ar@{.}[dr(0.8)]&&{}\save[]*{U_{\nu}}\restore\ar@{.}[ur(0.8)]\ar@{.}[2,2]&&\\
&&\bullet\ar@{.}[ur(0.8)]&&\ar[2,2]&&\bullet\ar@{.}[dr(0.8)]&\\
&&&&&\bullet\ar@{.}[dr(0.8)]\ar@{.}[ur(0.8)]&&\bullet\ar@{.}[dr(0.8)]&\\
&&&&\bullet\ar@{.}[dr(0.8)]\ar@{.}[ur(0.8)]&&{}\save[]*{U_1}\restore\ar@<-0.5ex>[2,2]\ar`dr/0pt[1,1][1,2]\ar@{.}[dr(0.8)]\ar@{.}[ur(0.8)]&&\bullet&\\
&&&&&\bullet\ar@{.}[r(0.8)]\ar@{.}[ur(0.8)]\ar@{.}[dr(0.8)]&\bullet\ar@{.}[r(0.6)]&\bullet\ar@{.}[dr(0.8)]\ar@{.}[ur(0.8)]&{}\save[]*{\quad W_+}\restore\ar@{.}[l(0.6)]\\
&&&&&&\bullet\ar@{.}[ur(0.8)]&&{}\save[]*{\quad W_-}\restore}
\end{array}
\]
}}where in the above picture the positions of $R$ and $W_{1}$ are fixed, however since  $U_{\nu}=V_{i_{\nu+1}}$ in general the length of the distance between $R$ and $U_{\nu}$ depends on $i_{\nu+1}$.  Also the distance between the other specials changes depending on $n-q$, however none of these changes in distance effect the irreducibility of the highlighted arrows. 

The above irreducible maps thus fix the positions of $W_1$, $W_+$, $W_-$ and the $U$'s, and the rest follow by noting the factorization of $(xy)^{i_{\nu+1}}$ as
\[
\xymatrix@C=65pt{ W_{i_{\nu+1}}&W_{i_{\nu+2}}\ar@/_0.65pc/[l]|{(xy)^{i_{\nu+1}-i_{\nu+2}}}\ar@{}[r]|{\hdots} & W_{i_{N-1}} & W_{1}\ar@/_0.65pc/[l]|{(xy)^{i_{N-1}-i_N}}&R\ar@/_0.65pc/[l]|{xy} }
\]
\end{proof}

The following is the main result of this section:
\begin{thm}\label{iso_of_rings2}
For a group $\mathbb{D}_{n,q}$ with parameter $0<\nu<N-1$, denote $R=\C{}[x,y]^{\mathbb{D}_{n,q}}$ and let $T_{n,q}=R\oplus W_{+}\oplus W_{-}\oplus_{s=1}^{\nu}U_s\oplus_{t=\nu+1}^{N}W_{i_t}$ be the sum of the special CM modules.  Then
\[
D_{n,q}\cong \t{End}_R(T_{n,q})\cong D_{n,q}^\prime.
\]
\end{thm}
\begin{proof}
We prove both statements at the same time by making different choices for the generators of the specials.   As stated in the proof of Lemma~\ref{specials2vertices2} the quiver of the endomorphism ring of the specials is precisely that of the quiver $Q$ defined above.  We first find representatives for the known number of arrows:

Set $m:=n-q$.  By choosing a $G$-equivariant basis as explained in Section 2, the following are always homomorphisms between the specials
{\tiny{
\[
\xymatrix@C=40pt@R=40pt{ &\bullet\ar@/^0.35pc/[d]^(0.4){\left(\begin{smallmatrix}
y^m&x^m \end{smallmatrix}\right)}&&&\\
\bullet\ar@/^0.35pc/[r]^(0.4){\left(\begin{smallmatrix}
\minus y^m&x^m \end{smallmatrix}\right)}&
\bullet\ar@/^0.35pc/[u]^(0.6){\left(\begin{smallmatrix}
x^m\\ \minus y^m \end{smallmatrix}\right)}\ar@/^0.35pc/[l]^(0.6){\left(\begin{smallmatrix}
x^m\\y^m \end{smallmatrix}\right)}\ar@/^0.35pc/[r]^{\left(\begin{smallmatrix}
\minus y^m&0\\0&x^m \end{smallmatrix}\right)}
&\bullet\ar@/^0.35pc/[l]^{\left(\begin{smallmatrix}
x^m&0\\0&y^m \end{smallmatrix}\right)}\ar@{}[r]|{\hdots} &\bullet\ar@/^0.35pc/[r]^{\left(\begin{smallmatrix}
\minus y^m&0\\0&x^m \end{smallmatrix}\right)}& \bullet\ar@/^0.35pc/[d]^(0.55){(xy)^{r_3}\left(\begin{smallmatrix}
\minus y^{i_{\nu+1}}\\ x^{i_{\nu+1}} \end{smallmatrix}\right)}\ar@/^0.35pc/[l]^{\left(\begin{smallmatrix}
x^m&0\\0&y^m \end{smallmatrix}\right)}\ar@/^0.35pc/[r]^{\left(\begin{smallmatrix}
y^{i_{\nu+1}}\\x^{i_{\nu+1}} \end{smallmatrix}\right)} & \bullet\ar@/^0.35pc/[l]^(0.2){(xy)^{r_3}\left(\begin{smallmatrix}
x^{i_{\nu+1}}&\minus y^{i_{\nu+1}} \end{smallmatrix}\right)}\ar@<-0.3ex>@{.}[r]&\bullet\ar@/_0.55pc/@{.}[l]\ar@{}[r]|{\hdots} & \bullet\ar@<-0.3ex>@{.}[r] & \bullet\ar@/_0.55pc/@{.}[l]\ar@<0.55ex>@{.}[1,-4] \\
&&&&\star\ar@/^0.35pc/[u]^(0.4){\left(\begin{smallmatrix}
x^{i_{\nu+1}}&y^{i_{\nu+1}} \end{smallmatrix}\right)}\ar@<-1.1ex>@{.}@/_0.55pc/[-1,4]&&}
\]
}}The proof that these maps are actually homomorphisms (equivalently $G$-equivariant) splits into two cases depending on the parity of $n-q$ since by definition both $G$ and the representations depend on $n-q$; we suppress the details.  By inspection they do not factor through any of the other specials and so can be chosen as representatives.  %\footnote{Also can view as compositions in the AR quiver, where its obvious that they don't factor through other specials.}
%\footnote{The signs on the $f$ arrows are ok - for $n-q$ odd just compose the preprojective choice made before, but for $n-q$ even verify directly using $G$-maps.  See ChoicesofArrows.tex for the proof.  Note this means that the composition in the twisted AR quiver (i.e. $n-q$ even) is slightly different, since otherwise we would get both signs positive on the $f$ arrows.}
Note that with these choices%\footnote{The $\pr{\nu+1}{\nu+1}$ are the same iff $\nu$ is even, similarly with the $\CLr{0}{\nu}k_1$.  The $\CLr{0}{\nu+1}$ are the same iff $\nu$ is odd.}
\begin{eqnarray*}
\pr{\nu+1}{\nu+1}&=&
\left\{
\begin{array}{l}
\ff{\nu+1}{\nu}\GG{\nu}{1}\ff{1}{+}\gg{+}{1}\FF{1}{\nu}\gg{\nu}{\nu+1}=(xy)^{r_3}w_2^{c_3}w_3^{d_3}=(xy)^{r_3}w_3\\
\frac{1}{2}\ff{\nu+1}{\nu}\GG{\nu}{1}(\ff{1}{+}\gg{+}{1}+\ff{1}{-}\gg{-}{1})\FF{1}{\nu}\gg{\nu}{\nu+1}=(xy)^{r_3}v_2^{c_3}v_3^{d_3}=(xy)^{r_3}v_3
\end{array}
\right.\\
\CLr{0}{\nu}k_1&=&
\left\{
\begin{array}{l}
\ff{0}{\nu}\GG{\nu}{1}\ff{1}{-}\gg{-}{1}\FF{1}{\nu}\gg{\nu}{0}=(xy)^{r_3}w_2^{c_3}w_3^{d_3}=(xy)^{r_3}w_3\\
\frac{1}{2}\ff{0}{\nu}\GG{\nu}{1}(\ff{1}{+}\gg{+}{1}+\ff{1}{-}\gg{-}{1})\FF{1}{\nu}\gg{\nu}{0}=(xy)^{r_3}v_2^{c_3}v_3^{d_3}=(xy)^{r_3}v_3
\end{array}
\right.\\
\CLr{0}{\nu+1}&=&
\left\{
\begin{array}{l}
 \ff{0}{\nu}\GG{\nu}{1}\ff{1}{+}\gg{+}{1}\FF{1}{\nu}\gg{\nu}{\nu+1}=w_2\\
\frac{1}{2}\ff{0}{\nu}\GG{\nu}{1}(\ff{1}{+}\gg{+}{1}+\ff{1}{-}\gg{-}{1})\FF{1}{\nu}\gg{\nu}{\nu+1}=v_2
\end{array}
\right.
\end{eqnarray*}
We must now choose representatives of the remaining arrows, and it is different choices of these which give the two different presentations.  As in the proof of the above lemma it is clear that the anti-clockwise arrows are $\an{N}{0}=xy$ and $\an{t}{t-1}=(xy)^{i_{t-1}-i_t}$.  Also, since we must reach the generators of the specials as paths from $\star$ (i.e. the vertex corresponding to $R$), by Lemma~\ref{1d_specials_generators} we can choose $\cl{t-1}{t}=w_2^{c_{l_{t-1}}}w_3^{d_{l_{t-1}}}$  for all $\nu+2\leq t\leq N$.  Alternatively, also by Lemma~\ref{1d_specials_generators}, we can choose $\cl{t-1}{t}=v_2^{c_{l_{t-1}}}v_3^{d_{l_{t-1}}}$  for all $\nu+2\leq t\leq N$ (these choices will give the symmetric presentation).  As explained in \cite[4.9]{Wemyss_reconstruct_D(i)} we may choose $\cl{N}{0}=w_2^{c_{l_N}}w_3^{d_{l_N}}=w_2^{c_{e-1}}w_3^{d_{e-1}}$ (respectively $v_2^{c_{e-1}}v_3^{d_{e-1}}$).

Thus we have proved that we may take as representatives of the arrows
{\tiny{
\[
\xymatrix@C=40pt@R=40pt{ &\bullet\ar@/^0.35pc/[d]^(0.4){\left(\begin{smallmatrix}
y^m&x^m \end{smallmatrix}\right)}&&&\\
\bullet\ar@/^0.35pc/[r]^(0.4){\left(\begin{smallmatrix}
\minus y^m&x^m \end{smallmatrix}\right)}&
\bullet\ar@/^0.35pc/[u]^(0.6){\left(\begin{smallmatrix}
x^m\\ \minus y^m \end{smallmatrix}\right)}\ar@/^0.35pc/[l]^(0.6){\left(\begin{smallmatrix}
x^m\\y^m \end{smallmatrix}\right)}\ar@/^0.35pc/[r]^{\left(\begin{smallmatrix}
\minus y^m&0\\0&x^m \end{smallmatrix}\right)}
&\bullet\ar@/^0.35pc/[l]^{\left(\begin{smallmatrix}
x^m&0\\0&y^m \end{smallmatrix}\right)}\ar@{}[r]|{\hdots} &\bullet\ar@/^0.35pc/[r]^{\left(\begin{smallmatrix}
\minus y^m&0\\0&x^m \end{smallmatrix}\right)}& \bullet\ar@/^0.35pc/[d]^(0.55){(xy)^{r_3}\left(\begin{smallmatrix}
\minus y^{i_{\nu+1}}\\ x^{i_{\nu+1}} \end{smallmatrix}\right)}\ar@/^0.35pc/[l]^{\left(\begin{smallmatrix}
x^m&0\\0&y^m \end{smallmatrix}\right)}\ar@/^0.35pc/[r]^{\left(\begin{smallmatrix}
y^{i_{\nu+1}}\\x^{i_{\nu+1}} \end{smallmatrix}\right)} & \bullet\ar@/^0.35pc/[l]^(0.2){(xy)^{r_3}\left(\begin{smallmatrix}
x^{i_{\nu+1}}&\minus y^{i_{\nu+1}} \end{smallmatrix}\right)}\ar@<-0.3ex>[r]|(0.45){\cl{\nu+1}{\nu+2}}&\bullet\ar@/_0.55pc/[l]_(0.4){(xy)^{i_{\nu+1}-i_{\nu+2}}}\ar@{}[r]|{\hdots} & \bullet\ar@<-0.3ex>[r]|{\cl{N-1}{N}} & \bullet\ar@/_0.55pc/[l]_(0.4){(xy)^{i_{N-1}-i_N}}\ar@<0.55ex>[1,-4]|{\cl{N}{0}} \\
&&&&\star\ar@/^0.35pc/[u]^(0.4){\left(\begin{smallmatrix}
x^{i_{\nu+1}}&y^{i_{\nu+1}} \end{smallmatrix}\right)}\ar@<-1.1ex>@/_0.55pc/[-1,4]_{xy}&&}
\]
}}with $\cl{t-1}{t}=w_2^{c_{l_{t-1}}}w_3^{d_{l_{t-1}}}$  for all $\nu+2\leq t\leq N$ and $\cl{N}{0}=w_2^{c_{e-1}}w_3^{d_{e-1}}$.   Now if $\sum_{i=1}^{N}(\alpha_i-2)\geq 2$ then for the moduli presentation we take %\footnote{$``k_1"=\ff{\nu+1}{\nu}\GG{\nu}{1}\ff{1}{-}\gg{-}{1}\FF{1}{\nu}\gg{\nu}{0} =(xy)^{2r_3}w_2^{c_2}w_3^{d_2}=(xy)^{2r_3}w_2$ (with obvious for symmetrical) is useful to get all invariants, since going anticlockwise $(xy)^{i_{\nu+1}}$ and using $r_2=2r_3+i_{\nu+1}$ (TF(v)) gives the invariant $(xy)^{r_2}w_2^{c_2}w_3^{d_2}$ at all vertices $\nu+1,\hdots,N$.  Note however going clockwise $w_2$ then $k_1$ isn't so sensible - get $(xy)^{2r_3}w_2^2$ which is only related to $(xy)^{2r_3}w_3^{2}$, which is an invariant squared.}
%\footnote{Also note the pattern that $\cl{N}{0}=k_{\sum(\alpha_i-2)+1}=k_{e-2}$, if $k_{\sum(\alpha_i-2)+1}$ was defined.}
\[
k_t=(xy)^{r_{t+2}}w_2^{c_{t+1}}w_3^{d_{t+1}} \mbox{ for all } 2\leq
t\leq \sum(\alpha_i-2)
\]
as representatives.  To see why we can do this, for example if $\alpha_{\nu+1}\geq 4$ consider $(xy)^{r_{4}}w_2^{c_{3}}w_3^{d_{3}}:W_{i_{\nu+1}}\rightarrow R$.  Firstly it does not factor through maps we have already chosen, since if it does then since $\ff{\nu+1}{\nu}$ has a $(xy)^{r_3}$ factor and $\cl{\nu+1}{\nu+2}=w_2^{c_{l_{\nu+1}}}w_3^{d_{l_{\nu+1}}}$ we may write $(xy)^{r_{4}}w_2^{c_{3}}w_3^{d_{3}}=(xy)^{r_3}F+w_2^{c_{l_{\nu+1}}}w_3^{d_{l_{\nu+1}}}h$ for some polynomials $F$ and $h$.  But by looking at $xy$ powers we know $(xy)^{r_4}$ divides $h$ and so after cancelling factors we may write $w_2^{c_3}w_3^{d_3}=(xy)^{r_3-r_4}F+w_2^{c_{l_{\nu+1}}}w_3^{d_{l_{\nu+1}}}h_1$.  After cancelling $w_2^{c_3}w_3^{d_3}$ (which $F$ must be divisible by) $1=(xy)^{r_3-r_4}F^\prime+w_2^{c_{l_{\nu+1}}-c_3}w_3^{d_{l_{\nu+1}}-d_3}h_1$ which is impossible since the assumption $\alpha_{\nu+1}\geq 4$ means that the right hand side cannot have degree zero terms.  %\footnote{since $\alpha_{\nu+1}>3$ implies that $l_{\nu+1}>3$.}
Secondly, $(xy)^{r_{4}}w_2^{c_{3}}w_3^{d_{3}}$ does not factor as a map $W_{i_{\nu+1}}\rightarrow R$ followed by a non-scalar invariant since this would contradict the embedding dimension.  The proof for the remainder of the arrows is similar; in fact the proof is identical to \cite[4.9]{Wemyss_reconstruct_D(i)}.

The symmetric presentation is the same, but we replace everywhere $w_2$ by $v_2$ and $w_3$ by $v_3$. 

The relations for the reconstruction algebra are now induced by the relations of matrix multiplication such that $x$ and $y$ commute. Using Lemma~\ref{r_is_difference_in_i_series}, Lemma~\ref{c_and_d_lemma} and the polynomial expressions for $\pr{\nu+1}{\nu+1}$, $\CLr{0}{\nu}k_1$, and $\CLr{0}{\nu+1}$ above, it is straightforward to verify that the relations in Definition~\ref{reconD2} (which we denote by $\c{S}^{\prime}$) are satisfied by these choices of polynomials, since the algorithmic relations are just the pattern in the cycles in $\tt{D}_1$ (respectively $\tt{D}_2$) from \cite[3.5]{Wemyss_reconstruct_D(i)}.

The remainder of the proof is now identical to \cite[4.9]{Wemyss_reconstruct_D(i)}: we first work in the completed case (so we can use \cite[3.3]{Wemyss_GL2} and \cite[3.4]{BIRS}) and we prove that the completion of the endomorphism ring of the specials is given as the completion of $\C{}Q$ (denoted $\C{}\hat{Q}$) modulo the closure of the ideal $\langle \c{S}^\prime\rangle$ (denoted $\overline{\langle \c{S}^\prime\rangle}$).  The non-completed version of the theorem then follows by simply taking the associated graded ring of both sides of the isomorphism.

Denote the kernel of the surjection $\C{}\hat{Q}\rightarrow \t{End}_{\C{}[[x,y]]^G}(T_{n,q}):=\Lambda$ by $I$, denote the radical of $\C{}\hat{Q}$ by $J$ and further for $t\in\{ \star,+,-,1,\hdots,N \}$ denote by $S_t$ the simple corresponding to the vertex $t$ of $Q$.  In Lemma~\ref{LIrelations2} below we show that the elements of $\c{S}^\prime$ are linearly independent in $I/ (IJ+JI)$.  Thus we may extend $\c{S}^\prime$ to a basis $\c{S}$ of $I/ (IJ+JI)$ and so by \cite[3.4(a)]{BIRS} $I=\overline{\langle \c{S}\rangle}$. But by combining \cite[3.4(b)]{BIRS}, \cite[3.3]{Wemyss_GL2} and inspecting our set $\c{S}^\prime$ we see that 
\[
\# (e_{a}\C{}\hat{Q}e_{b})\cap \c{S}=\# (e_{a}\C{}\hat{Q}e_{b})\cap \c{S}^\prime
\]
for all $a,b\in\{ \star,+,-,1,\hdots,N \}$, proving that the number of elements in $\c{S}$ and $\c{S}^\prime$ are the same.  Hence $\c{S}^\prime=\c{S}$ and so $I=\overline{\langle \c{S}^\prime\rangle}$, as required.
\end{proof}

The proof of the following lemma is only subtly different to that of \cite[4.10]{Wemyss_reconstruct_D(i)} and so we only prove the parts in which the subtle differences arise.
\begin{lemma}\label{LIrelations2}
With notation from the above proof, the members of $\c{S}^\prime$ are linearly independent in $I/(IJ+JI)$.
\end{lemma}
\begin{proof}
Following \cite[4.10]{Wemyss_reconstruct_D(i)} we say that a word $w$ in the path algebra $\C{}\hat{Q}$ satisfies condition (A) if
\[
\begin{array}{cl}
\t{(i)}&\t{It does not contain some proper subword which is a cycle.}\\
\t{(ii)}&\t{It does not contain some proper subword which is a path from $\star$ to $\nu+1$.}
\end{array}
\]
We know from the intersection theory and \cite[3.3]{Wemyss_GL2} that the ideal $I$ is generated by one relation from $\star$ to $\nu+1$, whereas all other generators are cycles.  Consequently if a word $w$ satisfies (A) then $w\notin IJ+JI$.  It is also clear that $\ff{0}{\nu}\gg{\nu}{\nu+1}-2\AN{0}{\nu+1}\notin IJ+JI$.

All members of $\c{S}^\prime$ are either cycles at some vertex or paths from $\star$ to $\nu+1$,  so to prove that the members of $\c{S}^\prime$ are linearly independent in $I/(IJ+JI)$ we just need to show that 
\begin{itemize}
\item[1.] the elements of $\c{S}^\prime$ which are paths from $\star$ to $\nu+1$ are linearly independent in $e_\star(I/(IJ+JI))e_{\nu+1}$.
\item[2.] for all $t\in\{ \star,+,-,1,\hdots,N \}$, the elements of $\c{S}^\prime$ that are cycles at $t$ are linearly independent in $e_t(I/(IJ+JI))e_t$.
\end{itemize}

The first condition is easy, since the only relation in $\c{S}^\prime$ from $\star$ to $\nu+1$ is $\ff{0}{\nu}\gg{\nu}{\nu+1}-2\AN{0}{\nu+1}$ and we have already noted that it does not belong to $IJ+JI$, thus it is non-zero and so linearly independent in $e_\star(I/(IJ+JI))e_{\nu+1}$.  Note also that there are only two paths of  minimal grade from $\star$ to $\nu+1$ (namely $\ff{0}{\nu}\gg{\nu}{\nu+1}$ and $\AN{0}{\nu+1}$) and by inspection of the polynomials they represent we do not have any other relation from $\star$ to $\nu+1$ of this grade.

For the second condition, we must check $t$ case by case: \\
\emph{Case $t=+$ and $t=-$}.  There is only one relation in $\c{S}^\prime$ which is a cycle at that vertex and so to prove linearly independence requires only that the relation is non zero in $I/(IJ+JI)$.  But clearly $\gg{+}{1}\ff{1}{+}\notin IJ+JI$ and $\gg{-}{1}\ff{1}{-}\notin IJ+JI$ since they both satisfy condition (A).\\
\emph{Case $t=1,\hdots\nu$}.  Again there is only one relation, and since each word satisfies condition (A) it does not belong to $IJ+JI$ and so the relation is non-zero.  For example if $\ff{1}{+}\gg{+}{1}-\ff{1}{-}\gg{-}{1}-2\ff{1}{2}\gg{2}{1}\in IJ+JI$ then $\ff{1}{+}\gg{+}{1}=\ff{1}{-}\gg{-}{1}+2\ff{1}{2}\gg{2}{1}+u$ for some $u\in IJ+JI$ in the free algebra $\C{}\hat{Q}$ which is impossible since $\ff{1}{+}\gg{+}{1}$ cannot appear in the right hand side.\\
\emph{Case $t=\nu+1$}.
If $\alpha_{\nu+1}=3$ then there are only two relations in $\c{S}^\prime$ from $\nu+1$ to $\nu+1$.  To prove linear independence suppose that
\[
\lambda_1(\cl{\nu+1}{\nu+2}\an{\nu+2}{\nu+1}-\pr{\nu+1}{\nu+1})+\lambda_2(\ff{\nu+1}{\nu}\gg{\nu}{\nu+1})=0
\]
in $e_{\nu+1}(I/(IJ+JI))e_{\nu+1}$.  Then
\[
\lambda_1\cl{\nu+1}{\nu+2}\an{\nu+2}{\nu+1}+\lambda_2\ff{\nu+1}{\nu}\gg{\nu}{\nu+1}=\lambda_1\pr{\nu+1}{\nu+1}+u
\]
in the free algebra $\C{}\hat{Q}$ for some $u\in IJ+JI$.  But $\cl{\nu+1}{\nu+2}\an{\nu+2}{\nu+1}$ satisfies condition (A) and so cannot appear on the right hand side, forcing $\lambda_1=0$. Similarly  since $\ff{\nu+1}{\nu}\gg{\nu}{\nu+1}$ satisfies condition (A) $\lambda_2=0$.  This proves the case $\alpha_{\nu+1}=3$ and so we can assume that $\alpha_{\nu+1}>3$.  Now there are $\alpha_{\nu+1}-1$ relations in $\c{S}^\prime$ from $\nu+1$ to $\nu+1$, and suppose that%\footnote{If $\alpha_{\nu+1}>3$ then $u_{\nu+1}=\alpha_{\nu+1}-2$ since we start labelling at $k_2$.}
\begin{multline*}
\lambda_1(\ff{\nu+1}{\nu}\gg{\nu}{\nu+1})+\lambda_2(k_2\AN{0}{\nu+1}-\pr{\nu+1}{\nu+1})\\ +\sum_{p=2}^{u_{\nu+1}-1}\lambda_{p+1}(k_{p+1}\AN{0}{\nu+1}-k_p\CLr{0}{\nu+1})+\lambda_{u_{\nu+1}+1}(\cl{\nu+1}{\nu+2}\an{\nu+2}{\nu+1}-k_{u_{\nu+1}}\CLr{0}{\nu+1})=0
\end{multline*}
in $e_{\nu+1}(I/(IJ+JI))e_{\nu+1}$ (where the sum may be empty).  Then%\footnote{The $k_2$ is absorbed into the left hand sum}
\begin{multline*}
\lambda_1\ff{\nu+1}{\nu}\gg{\nu}{\nu+1}+\sum_{p=1}^{u_{\nu+1}-1}\lambda_{p+1}k_{p+1}\AN{0}{\nu+1}+\lambda_{u_{\nu+1}+1}\cl{\nu+1}{\nu+2}\an{\nu+2}{\nu+1}\\=\lambda_2\pr{\nu+1}{\nu+1} +\sum_{p=2}^{u_{\nu+1}-1}\lambda_{p+1}k_p\CLr{0}{\nu+1}+\lambda_{u_{\nu+1}+1}k_{u_{\nu+1}}\CLr{0}{\nu+1}+u
\end{multline*}
in the free algebra $\C{}\hat{Q}$ for some $u\in IJ+JI$.  Now $\ff{\nu+1}{\nu}\gg{\nu}{\nu+1}$ and $\cl{\nu+1}{\nu+2}\an{\nu+2}{\nu+1}$ satisfy condition (A) forcing  $\lambda_1=0$ and $\lambda_{u_{\nu+1}+1}=0$.  Thus
\[
\lambda_{u_{\nu+1}}k_{u_{\nu+1}}\AN{0}{\nu+1}\equiv \lambda_2\pr{\nu+1}{\nu+1}+\t{terms starting with k of strictly smaller index}
\]
mod $IJ+JI$.  But since $k_{u_{\nu+1}}\AN{0}{\nu+1}$ does not have any subwords which are cycles, the only way we can change it mod $IJ+JI$ is to bracket as $k_{u_{\nu+1}}(\AN{0}{\nu+1})$ and use the relation in $I$ from $\star$ to $\nu+1$.  %\footnote{The $\AN{}{}$'s are the main difference between the $\nu=0$ case and here - we have a relation which messes up the $\AN{}{}$'s in a different way.  The biggest difference is $\pr{\nu+1}{\nu+1}$ though, since it contains a subword which is a cycle.}
Doing this we get $k_{u_{\nu+1}}\AN{0}{\nu+1}\equiv \frac{1}{2}k_{u_{\nu+1}}\ff{0}{\nu}\gg{\nu}{\nu+1}$.  This still does not start with $\pr{\nu+1}{\nu+1}$ or $k$ of strictly lower index, so we must again use relations in $IJ+JI$ to change the terms.  But $k_{u_{\nu+1}}\ff{0}{\nu}\gg{\nu}{\nu+1}$ does not contain any subwords which are cycles, which means the only way to change it mod $IJ+JI$ is to use a relation for $\ff{0}{\nu}\gg{\nu}{\nu+1}$.  But there is only one relation for $\ff{0}{\nu}\gg{\nu}{\nu+1}$ in $I$ and so using it we arrive back at
\[
k_{u_{\nu+1}}\AN{0}{\nu+1}\equiv \frac{1}{2}k_{u_{\nu+1}}\ff{0}{\nu}\gg{\nu}{\nu+1}\equiv k_{u_{\nu+1}}\AN{0}{\nu+1}
\]
mod $IJ+JI$.  Thus mod $IJ+JI$ it is impossible to transform $k_{u_{\nu+1}}\AN{0}{\nu+1}$ into an expression involving $\pr{\nu+1}{\nu+1}$ or $k$ terms with strictly lower index, thus we must have $\lambda_{u_{\nu+1}}=0$.   Now re-arranging we get
\[
\lambda_{u_{\nu+1}-1}k_{u_{\nu+1}-1}\AN{0}{\nu+1}\equiv \lambda_2\pr{\nu+1}{\nu+1}+\t{terms starting with k of strictly smaller index}
\]
mod $IJ+JI$ and so by induction all the $\lambda$'s are zero, as required.\\
\emph{Case $t$ with $\nu+2\leq t\leq N$.}  %\footnote{again, as in $\mu=0$ can get by without doing $N$ separately} 
The proof here is identical to the $2\leq t\leq N$ case in the proof of \cite[4.10]{Wemyss_reconstruct_D(i)}.  %\footnote{Since the $A$'s now satisfy condition (A).}
\\
\emph{Case $t=\star$}. This is proved using a very similar argument as in Case $t=\nu+1$ above.
\end{proof}

\section{The reconstruction algebra for $\nu=N-1$}
In this section we define the reconstruction algebra $D_{n,q}$ when $\nu=N-1$ and prove that it is isomorphic to the endomorphism ring of the special CM modules.  This is in fact very similar to the previous section, but here the vertex $\star$ is connected in a slightly different way and also the number of extra arrows out of the vertex $\nu+1$ is $\alpha_{\nu+1}-2$ compared to $\alpha_{\nu+1}-3$ previously. 

Consider, for $N\in\mathbb{N}$ with $N\geq 2$ and for a positive integer $\alpha_N\geq 2$, the labelled Dynkin diagram of type D:
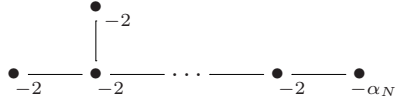
\[
\xymatrix@C=20pt@R=15pt{ &\bullet\ar@{-}[d]^<{-2}&&&\\
\bullet\ar@{-}[r]_<{-2} & \bullet\ar@{-}[r]_<{-2}
&\hdots\ar@{-}[r] &\bullet\ar@{-}[r]_<{-2} & \bullet
\ar@{}[r]_<{-\alpha_N}&}
\]
We call the left hand vertex the $+$ vertex, the top vertex the $-$ vertex and the remaining vertices $1,\hdots,N$ reading from left to right. To this picture we add an extended vertex $\star$ and
`double-up' as follows:
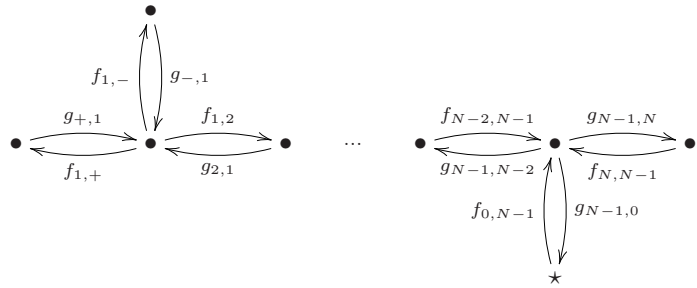
\[
\xymatrix@C=40pt@R=40pt{ &\bullet\ar@/^0.35pc/[d]^{\gg{-}{1}}&&&\\
\bullet\ar@/^0.35pc/[r]^{\gg{+}{1}}&
\bullet\ar@/^0.35pc/[u]^{\ff{1}{-}}\ar@/^0.35pc/[l]^{\ff{1}{+}}\ar@/^0.35pc/[r]^{\ff{1}{2}}
&\bullet\ar@/^0.35pc/[l]^{\gg{2}{1}}\ar@{}[r]|{\hdots} &\bullet\ar@/^0.35pc/[r]^{\ff{N-2}{N-1}}& \bullet\ar@/^0.35pc/[l]^{\gg{N-1}{N-2}}\ar@/^0.35pc/[d]^{\gg{N-1}{0}}\ar@/^0.35pc/[r]^{\gg{N-1}{N}} & \bullet\ar@/^0.35pc/[l]^{\ff{N}{N-1}} \\
&&&&\star\ar@/^0.35pc/[u]^{\ff{0}{N-1}}&}
\]
This is precisely the underlying quiver of the preprojective algebra of type $D$.  Now if $\alpha_N>2$ then add an extra $\alpha_N-2$ arrows from the $N^{th}$ vertex to $\star$ and label them $k_1,\hdots,k_{\alpha_N-2}$ reading from
left to right.  Call this quiver $Q$.
\begin{example}\label{D74quiver}
\t{Consider $\mathbb{D}_{7,4}$ then $\frac{7}{4}=[2,4]$ and
so $Q$ is}
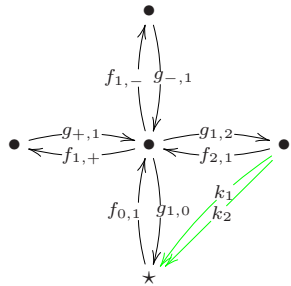
\[
\xymatrix@C=40pt@R=40pt{ &\bullet\ar@/^0.35pc/[d]|{\quad\gg{-}{1}}&&\\
\bullet\ar@/^0.35pc/[r]|{\gg{+}{1}}&
\bullet\ar@/^0.35pc/[d]|{\quad\gg{1}{0}}\ar@/^0.35pc/[u]|{\ff{1}{-}\quad}\ar@/^0.35pc/[l]|{\ff{1}{+}}\ar@/^0.35pc/[r]|{\gg{1}{2}}&
 \bullet\ar@/^0.35pc/[l]|{\ff{2}{1}}\ar@[green]@<0.25ex>@/_0.35pc/[1,-1]|(0.4){k_1}\ar@[green]@<0.25ex>[1,-1]|{\,\,k_2} \\
&\star\ar@/^0.35pc/[u]|{\ff{0}{1}\quad}&&}
\]
\end{example}
\begin{example}\label{D75quiver}
\t{Consider $\mathbb{D}_{7,5}$ then $\frac{7}{5}=[2,2,3]$ and so $Q$ is}
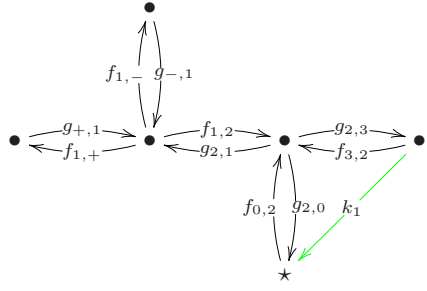
\[
\xymatrix@C=40pt@R=40pt{ &\bullet\ar@/^0.35pc/[d]|{\quad\gg{-}{1}}&&\\
\bullet\ar@/^0.35pc/[r]|{\gg{+}{1}}&
\bullet\ar@/^0.35pc/[u]|{\ff{1}{-}\quad}\ar@/^0.35pc/[l]|{\ff{1}{+}}\ar@/^0.35pc/[r]|{\ff{1}{2}}
&\bullet\ar@/^0.35pc/[d]|{\quad\gg{2}{0}}\ar@/^0.35pc/[l]|{\gg{2}{1}}\ar@/^0.35pc/[r]|{\gg{2}{3}} & \bullet\ar@/^0.35pc/[l]|{\ff{3}{2}}\ar@[green][1,-1]|{k_1} \\
&&\star\ar@/^0.35pc/[u]|{\ff{0}{2}\quad}&&}
\]
\end{example}

\begin{defin}\label{reconD1}
For $\frac{n}{q}=[\alpha_1,\hdots,\alpha_N]$ with $\nu=N-1$ define
$D_{n,q}$, the reconstruction algebra of type $D$, to be the path
algebra of the quiver $Q$ defined above subject to the relations
\[
\begin{array}{rl}
\mbox{If }N=2& \left\{\begin{array}{l}\gg{+}{1}\ff{1}{+}=0\quad \ff{0}{1}\gg{1}{0}=0 \\
 \gg{-}{1}\ff{1}{-}=0 \quad \ff{2}{1}\gg{1}{2}=0 \\
\gg{1}{0}\ff{0}{1}-\gg{1}{2}\ff{2}{1}=\ff{1}{+}\gg{+}{1}-\ff{1}{-}\gg{-}{1}\end{array}\right.\\
&\\
\mbox{If }N>2& \left\{\begin{array}{l}\gg{+}{1}\ff{1}{+}=0\quad \ff{0}{N-1}\gg{N-1}{0}=0\\  
\gg{-}{1}\ff{1}{-}=0\quad \ff{N}{N-1}\gg{N-1}{N}=0 \\
\ff{1}{+}\gg{+}{1}-\ff{1}{-}\gg{-}{1}=2\ff{1}{2}\gg{2}{1}\\
\gg{t}{t-1}\ff{t-1}{t}=\ff{t}{t+1}\gg{t+1}{t}\mbox{ for all }2\leq t\leq N-2\\
\gg{N-1}{0}\ff{0}{N-1}-\gg{N-1}{N}\ff{N}{N-1}=2\gg{N-1}{N-2}\ff{N-2}{N-1}\end{array}\right.
\end{array}
\]
Further if $\alpha_N>2$ we add, in each case, the relations 
\[
\begin{array}{l}
k_1\AN{0}{N}=\pr{N}{N}, \AN{0}{N}k_1=\pr{0}{0}\\
k_t\CLr{0}{N}=k_{t+1}\AN{0}{N}, \AN{0}{N}k_{t+1}=\CLr{0}{N}k_t \mbox{ for all }1\leq t\leq \alpha_N-3
\end{array}
\]
where $\AN{0}{N}:=\frac{1}{2}\ff{0}{N-1}\gg{N-1}{N}$.  The only thing that remains to be defined is the $\pr{}{}$'s and $\CLr{}{}$'s, and it is these which change according to the presentation, namely%\footnote{$\pr{N}{N}$ represents same things as in previous footnote.  Don't have all the $k_1$ or $\CLr{0}{N-1}=\CLr{0}{\nu}$ nonsense here since just explicitly state it in the relations as $\pr{0}{0}$, since if it happens it always happens at vertex $N$! (c.f. before, problem was going back when at different vertex). $\CLr{0}{N}$ represents same things as before}
\begin{eqnarray*}
\pr{0}{0}&=&\left\{  \begin{array}{c} \ff{0}{N-1}\GG{N-1}{1}\ff{1}{-}\gg{-}{1}\FF{1}{N-1}\gg{N-1}{0} \\ \frac{1}{2}\ff{0}{N-1}\GG{N-1}{1}(\ff{1}{+}\gg{+}{1}+\ff{1}{-}\gg{-}{1})\FF{1}{N-1}\gg{N-1}{0}\end{array} \right.\\ 
\pr{N}{N}&=&\left\{  \begin{array}{c} \ff{N}{N-1}\GG{N-1}{1}\ff{1}{+}\gg{+}{1}\FF{1}{N-1}\gg{N-1}{N} \\ \frac{1}{2}\ff{N}{N-1}\GG{N-1}{1}(\ff{1}{+}\gg{+}{1}+\ff{1}{-}\gg{-}{1})\FF{1}{N-1}\gg{N-1}{N} \end{array} \right.\\ 
\CLr{0}{N}&=&\left\{  \begin{array}{c}  \ff{0}{N-1}\GG{N-1}{1}\ff{1}{+}\gg{+}{1}\FF{1}{N-1}\gg{N-1}{N} \\ \frac{1}{2}\ff{0}{N-1}\GG{N-1}{1}(\ff{1}{+}\gg{+}{1}+\ff{1}{-}\gg{-}{1})\FF{1}{N-1}\gg{N-1}{N} \end{array} \right.
\end{eqnarray*}
where for the moduli presentation we take the top choices and for the symmetric presentation we take the bottom choices.
\end{defin}
\begin{remark}
\t{The relations are very similar to before, but now all relations are cycles.  As in the previous section for the moment denote $D_{n,q}$ for the moduli presentation and $D_{n,q}^\prime$ for the symmetric presentation; we shall see in Theorem~\ref{iso_of_rings1} that they are isomorphic.}%\footnote{The presentations are exactly the same if and only if either (i) in base case SL or (ii) We have $\alpha_N=3$ with $\nu$ even (since then $\pr{N}{N}$ and $\pr{0}{0}$ are the same in both presentations, and $\CLr{0}{N}$ doesn't enter the picture)}
\end{remark}
\begin{example}\label{D74relations}
\t{For the group $\mathbb{D}_{7,4}$ the moduli presentation of the reconstruction algebra is the quiver in Example~\ref{D74quiver} subject to the relations}
\[
\begin{array}{c}
\begin{array}{l}
 \gg{+}{1}\ff{1}{+}=0\\ \gg{-}{1}\ff{1}{-}=0 \\ \ff{0}{1}\gg{1}{0}=0 \\ \ff{2}{1}\gg{1}{2}=0 \end{array}\\
 \gg{1}{0}\ff{0}{1}-\gg{1}{2}\ff{2}{1}=\ff{1}{+}\gg{+}{1}-\ff{1}{-}\gg{-}{1}\\
\begin{array}{cc}
\begin{array}{r}
\ff{2}{1}\ff{1}{+}\gg{+}{1}\gg{1}{2}=k_1(\frac{1}{2}\ff{0}{1}\gg{1}{2})\\
k_1(\ff{0}{1}\ff{1}{+}\gg{+}{1}\gg{1}{2})=k_2(\frac{1}{2}\ff{0}{1}\gg{1}{2})
\end{array}&
\begin{array}{l}
(\frac{1}{2}\ff{0}{1}\gg{1}{2})k_1=\ff{0}{1}\ff{1}{-}\gg{-}{1}\gg{1}{0}\\
(\frac{1}{2}\ff{0}{1}\gg{1}{2})k_2=(\ff{0}{1}\ff{1}{+}\gg{+}{1}\gg{1}{2})k_1
\end{array}
\end{array}
\end{array}
\]
\end{example}
\begin{example}\label{D75relations}
\t{For the group $\mathbb{D}_{7,5}$ the symmetric presentation of the reconstruction algebra is the quiver in Example~\ref{D75quiver} subject to the relations}
\[
\begin{array}{cc}
\begin{array}{l}
\gg{+}{1}\ff{1}{+}=0\\ \gg{-}{1}\ff{1}{-}=0\\ \ff{0}{2}\gg{2}{0}=0\\ \ff{3}{2}\gg{2}{3}=0 \end{array}\\
\begin{array}{c}\ff{1}{+}\gg{+}{1}-\ff{1}{-}\gg{-}{1}=2\ff{1}{2}\gg{2}{1}\\
\gg{2}{0}\ff{0}{2}-\gg{2}{3}\ff{3}{2}=2\gg{2}{1}\ff{1}{2}
\end{array}\\
\begin{array}{c}
\frac{1}{2}\ff{3}{2}\gg{2}{1}(\ff{1}{+}\gg{+}{1}+\ff{1}{-}\gg{-}{1})\ff{1}{2}\gg{2}{3}=k_1(\frac{1}{2}\ff{0}{2}\gg{2}{3})\\
 \frac{1}{2}\ff{0}{2}\gg{2}{1}(\ff{1}{+}\gg{+}{1}+\ff{1}{-}\gg{-}{1})\ff{1}{2}\gg{2}{0}=(\frac{1}{2}\ff{0}{2}\gg{2}{3})k_1\end{array}
\end{array}
\]
\end{example}

We now justify the positions of the specials relative to the dual graph.  The proof is very similar to that of Lemma~\ref{specials2vertices2} but is sufficiently illuminating to be worthy of inclusion.
\begin{lemma}\label{specials2vertices1}
Consider $\mathbb{D}_{n,q}$ with $\nu=N-1$.  Then the specials correspond to the dual graph of the minimal resolution in the following way
\[
\begin{array}{cc}
\begin{array}{c}
\xymatrix@C=20pt@R=15pt{ &\bullet\ar@{-}[d]^<{-2}&&&\\
\bullet\ar@{-}[r]_<{-2} & \bullet\ar@{-}[r]_<{-2}
&\hdots\ar@{-}[r] &\bullet\ar@{-}[r]_<{-2} & \bullet
\ar@{}[r]_<{-\alpha_N}&}
\end{array} &
\begin{array}{c}
\xymatrix@C=15pt@R=10pt{ &W_-\ar@{-}[d]&&&\\
W_+\ar@{-}[r] & U_{1}\ar@{-}[r]
&\hdots\ar@{-}[r] &U_{N-1}\ar@{-}[r] & W_{i_N}}
\end{array}
\end{array}
\]
\end{lemma}
\begin{proof}
Now $Z_f=E_{+}+E_{-}+\sum_{p=1}^{N-1}2E_p+E_N$ as in the $\t{SL}(2,\C{})$ case.  Thus by \cite[3.6]{Wemyss_GL2} the quiver of the endomorphism ring of the specials is precisely the $Q$ defined above in this section. 

It is clear that the following composition of maps in the AR quiver do not factor through other specials
{\tiny{
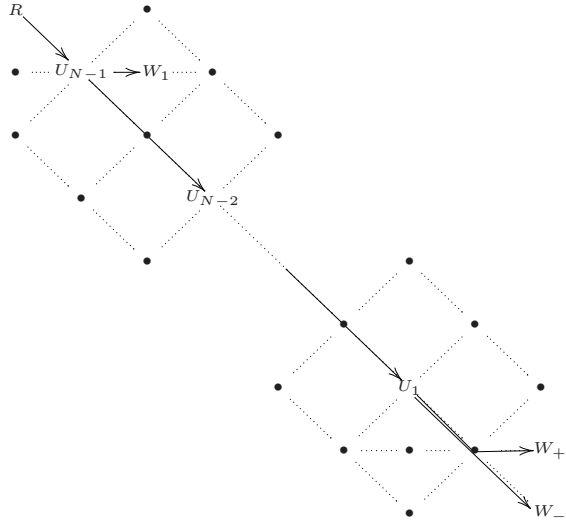
\[
\begin{array}{c}
\xymatrix@R=15pt@C=15pt{{}\save[]*{R}\restore\ar[dr(0.8)] &&\bullet\ar@{.}[dr(0.8)]&\\
\bullet\ar@{.}[r(0.5)]&{}\save[]*{U_{N-1}}\restore\ar[2,2]\ar@{.}[ur(0.8)]\ar@{.}[dr(0.8)]&{}\save[]*{\quad W_1}\restore\ar@{<-}[l(0.5)]&\bullet\ar@{.}[l(0.6)]\ar@{.}[dr(0.8)]&\\
\bullet\ar@{.}[ur(0.8)]\ar@{.}[dr(0.8)]&&\bullet\ar@{.}[ur(0.8)]\ar@{.}[dr(0.8)]&&\bullet&&\\
&\bullet\ar@{.}[ur(0.8)]\ar@{.}[dr(0.8)]&&{}\save[]*{U_{N-2}}\restore\ar@{.}[ur(0.8)]\ar@{.}[2,2]&&\\
&&\bullet\ar@{.}[ur(0.8)]&&\ar[2,2]&&\bullet\ar@{.}[dr(0.8)]&\\
&&&&&\bullet\ar@{.}[dr(0.8)]\ar@{.}[ur(0.8)]&&\bullet\ar@{.}[dr(0.8)]&\\
&&&&\bullet\ar@{.}[dr(0.8)]\ar@{.}[ur(0.8)]&&{}\save[]*{U_1}\restore\ar@<-0.5ex>[2,2]\ar`dr/0pt[1,1][1,2]\ar@{.}[dr(0.8)]\ar@{.}[ur(0.8)]&&\bullet&\\
&&&&&\bullet\ar@{.}[r(0.8)]\ar@{.}[ur(0.8)]\ar@{.}[dr(0.8)]&\bullet\ar@{.}[r(0.6)]&\bullet\ar@{.}[dr(0.8)]\ar@{.}[ur(0.8)]&{}\save[]*{\quad W_+}\restore\ar@{.}[l(0.6)]\\
&&&&&&\bullet\ar@{.}[ur(0.8)]&&{}\save[]*{\quad W_-}\restore}
\end{array}
\]
}}where the difference now (compared to Lemma~\ref{specials2vertices2}) is that in the above picture the positions of $R$, $W_{1}$ \emph{and} $U_{N-1}$ are fixed; again the length of the distance between the positions of the other specials changes depending on $n-q$ (in the picture we've drawn $n-q=2$) but this does not effect the fact that the maps are irreducible.  
\end{proof}
\begin{remark}
\t{The picture in the above proof explains algebraically why in the case $\nu=N-1$ we don't have to directly connect $R$ to the vertex $W_{i_N}=W_1$, since any such map must necessarily factor through $U_{N-1}$.}
\end{remark}

The following result is the main result of this section:
\begin{thm}\label{iso_of_rings1}
For a group $\mathbb{D}_{n,q}$ with parameter $\nu=N-1$, denote $R=\C{}[x,y]^{\mathbb{D}_{n,q}}$ and let $T_{n,q}=R\oplus W_{+}\oplus W_{-}\oplus_{s=1}^{N-1}U_s\oplus W_{i_N}$ be the sum of the special CM modules.  Then
\[
D_{n,q}\cong \t{End}_R(T_{n,q})\cong D_{n,q}^\prime.
\]
\end{thm}
\begin{proof}
The proof is very similar to that of Theorem~\ref{iso_of_rings2}.  Setting $m:=n-q$ then as before we may choose the following as representatives of the irreducible maps between the specials %\footnote{The preprojective part is exactly the same as before, since here $i_{\nu+1}=1$}
{\tiny{
\[
\xymatrix@C=40pt@R=40pt{ &\bullet\ar@/^0.35pc/[d]^(0.4){\left(\begin{smallmatrix}
y^m&x^m \end{smallmatrix}\right)}&&&\\
\bullet\ar@/^0.35pc/[r]^(0.4){\left(\begin{smallmatrix}
\minus y^m&x^m \end{smallmatrix}\right)}&
\bullet\ar@/^0.35pc/[u]^(0.6){\left(\begin{smallmatrix}
x^m\\ \minus y^m \end{smallmatrix}\right)}\ar@/^0.35pc/[l]^(0.6){\left(\begin{smallmatrix}
x^m\\y^m \end{smallmatrix}\right)}\ar@/^0.35pc/[r]^{\left(\begin{smallmatrix}
\minus y^m&0\\0&x^m \end{smallmatrix}\right)}
&\bullet\ar@/^0.35pc/[l]^{\left(\begin{smallmatrix}
x^m&0\\0&y^m \end{smallmatrix}\right)}\ar@{}[r]|{\hdots} &\bullet\ar@/^0.35pc/[r]^{\left(\begin{smallmatrix}
\minus y^m&0\\0&x^m \end{smallmatrix}\right)}& \bullet\ar@/^0.35pc/[l]^{\left(\begin{smallmatrix}
x^m&0\\0&y^m \end{smallmatrix}\right)}\ar@/^0.35pc/[d]^(0.6){(xy)^{r_3}\left(\begin{smallmatrix}
\minus y\\ x \end{smallmatrix}\right)}\ar@/^0.35pc/[r]^{\left(\begin{smallmatrix}
y\\x \end{smallmatrix}\right)} & \bullet\ar@/^0.35pc/[l]^(0.4){(xy)^{r_3}\left(\begin{smallmatrix}
x&\minus y \end{smallmatrix}\right)}\\
&&&&\star\ar@/^0.35pc/[u]^(0.4){\left(\begin{smallmatrix}
x&y\end{smallmatrix}\right)}&}
\]
}}We note that with these choices we have
\begin{eqnarray*}
\pr{N}{N}&=&
\left\{
\begin{array}{l}
\ff{N}{N-1}\GG{N-1}{1}\ff{1}{+}\gg{+}{1}\FF{1}{N-1}\gg{N-1}{N}=(xy)^{r_3}w_2^{c_3}w_3^{d_3}=(xy)^{r_3}w_3\\
\frac{1}{2}\ff{N}{N-1}\GG{N-1}{1}(\ff{1}{+}\gg{+}{1}+\ff{1}{-}\gg{-}{1})\FF{1}{N-1}\gg{N-1}{N}=(xy)^{r_3}v_2^{c_3}v_3^{d_3}=(xy)^{r_3}v_3
\end{array}
\right.\\
\pr{0}{0}&=&
\left\{
\begin{array}{l}
\ff{0}{N-1}\GG{N-1}{1}\ff{1}{-}\gg{-}{1}\FF{1}{N-1}\gg{N-1}{0}=(xy)^{r_3}w_2^{c_3}w_3^{d_3}=(xy)^{r_3}w_3\\
\frac{1}{2}\ff{0}{N-1}\GG{N-1}{1}(\ff{1}{+}\gg{+}{1}+\ff{1}{-}\gg{-}{1})\FF{1}{N-1}\gg{N-1}{0}=(xy)^{r_3}v_2^{c_3}v_3^{d_3}=(xy)^{r_3}v_3
\end{array}
\right.\\
\CLr{0}{N}&=&
\left\{
\begin{array}{l}
 \ff{0}{N-1}\GG{N-1}{1}\ff{1}{+}\gg{+}{1}\FF{1}{N-1}\gg{N-1}{N}=w_2\\
\frac{1}{2}\ff{0}{N-1}\GG{N-1}{1}(\ff{1}{+}\gg{+}{1}+\ff{1}{-}\gg{-}{1})\FF{1}{N-1}\gg{N-1}{N}=v_2
\end{array}
\right. 
\end{eqnarray*}
and $\AN{0}{N}=xy$.  Now if $\alpha_N>2$ there are extra $k$ arrows, and here in the moduli presentation we choose for representatives %\footnote{adding 3 and 2 here, unlike the other two $\nu$ cases!}
\[
k_t= (xy)^{r_{t+3}}w_2^{c_{t+2}}w_3^{d_{t+2}}\mbox{ for all } 1\leq
t\leq \alpha_N-2.
\]
The proof that we can actually choose these is identical to before.  %\footnote{Actually its easier, since now the only map to other specials out of $N$ has $(xy)^{r_3}$ factor, so don't have the $\cl{\nu+1}{\nu+2}$ arrow in the way.}
Again the symmetric presentation is obtained by everywhere replacing $w_2$ by $v_2$ and $w_3$ by $v_3$.

The relations for the reconstruction algebra are again induced by matrix multiplication such that $x$ and $y$ commute. By Lemma~\ref{r_is_difference_in_i_series} and Lemma~\ref{c_and_d_lemma} it is very easy to verify that the algorithmic relations in Definition~\ref{reconD1} are satisfied with these choices of polynomials.  From here the proof is identical to that of Theorem~\ref{iso_of_rings2} (actually all relations are now cycles, so in fact the proof of linear independence becomes a little easier) hence we suppress the details.
\end{proof}

\section{Moduli Examples}
In this section we justify some of the philosophy given in the introduction and also in the introduction of \cite{Wemyss_reconstruct_D(i)}.  Geometrically the key change in viewpoint is that we should not view the minimal resolution as $G$-Hilb (which we can do via a result of Ishii \cite{Ishii}), but rather we should instead view the minimal resolution as being very similar to a space that we already understand.  It is the reconstruction algebra which tells us which space to compare to, and it is the reconstruction algebra which encodes the difference. 

In this section (to ease notation) we shall use reconstruction algebras to compare the explicit geometry of the following two examples
\[
\begin{array}{cc}
\begin{array}{c}
\xymatrix@C=15pt@R=10pt{ &\bullet\ar@{-}[d]^<{\minus 2}&&\\
\bullet\ar@{-}[r]_<{\minus 2} & \bullet\ar@{-}[r]_<{\minus 2}
& \bullet\ar@{}[r]_<{\minus 2}&}\end{array}&\begin{array}{c}
\xymatrix@C=15pt@R=10pt{ &\bullet\ar@{-}[d]^<{\minus 2}&&\\
\bullet\ar@{-}[r]_<{\minus 2} & \bullet\ar@{-}[r]_<{\minus 2}
& \bullet\ar@{}[r]_<{\minus 3}&}\end{array}\end{array}
\]
from which the pattern is clear.  The left example corresponds to the group $\mathbb{D}_{3,2}$ whereas the right example corresponds to $\mathbb{D}_{5,2}$. The first is an example inside $SL(2,\C{})$ and so has been extensively studied by many people.  Below we will produce an open cover of the minimal resolution of  $\C{2}/\mathbb{D}_{3,2}$ which we draw as follows:
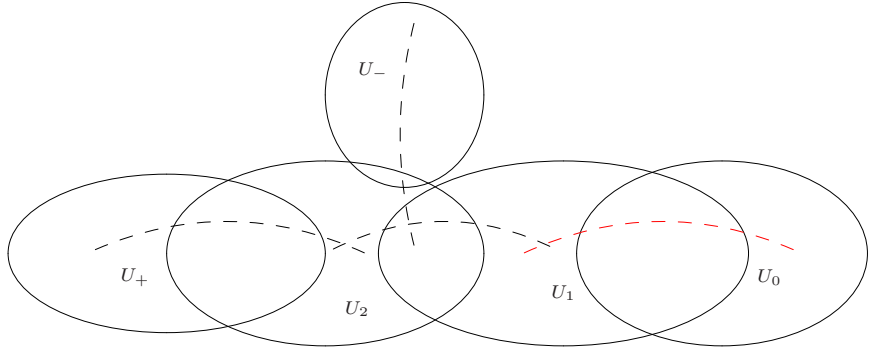
\[
\xymatrix{&&&&\ar@/_0.5pc/@<1ex>@{--}[3,0]&&&&&\\
&&&\ar@{}[-1,2]|(0.3){U_-}&{}\drop\xycircle<30pt,35pt>{}&&&&\\&&&&&&&\\
\ar@{}[1,2]|(0.3){U_+}{}\ar@{--}@/^1pc/[0,4]+<-15pt,0pt>&{}\drop\xycircle<60pt,30pt>{}&\ar@{}[1,2]|(0.7){U_2}&
\ar@{--}@/^1pc/[0,3] {}\drop\xycircle<60pt,35pt>{}&&\ar@{}[1,2]|(0.5){U_1}&{}\drop\xycircle<70pt,35pt>{}&
&\ar@{}[1,2]|(0.3){U_0}{}\drop\xycircle<55pt,35pt>{}&{}\ar@[red]@{--}@/_1pc/[0,-4]+<15pt,0pt>\\&&&&&&&&&&&}
\]
The point is that by changing reconstruction algebra we change the red (-2)-curve in the minimal resolution into a $(\minus 3)$ curve in a very efficient way.  Notice first that the fundamental cycles in the two examples are the same and consequently the reconstruction algebras are very similar.  We shall see that the reconstruction algebra \emph{only} changes the equation of the open set $U_0$ and the glue between $U_0$ and $U_1$.  The change in glue gives the change in self-intersection number.  The rest of the open sets (and their glues) will remain the same and so we will have the desired configuration of $\mathbb{P}^1$s.  This is geometrically the nicest solution since we are changing the least ammount of information to obtain one space from the other.  Note that if we use the $G$-Hilb description of the minimal resolution (or the McKay quiver) it is not obvious that this should be true.

\begin{example}
\t{Consider the group $\mathbb{D}_{3,2}$ of order 8.  Since this is inside $\t{SL}(2,\C{})$ the reconstruction algebra is just the preprojective algebra
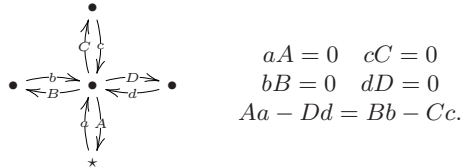
\[
\begin{array}{cc}
{\scriptsize{\begin{array}{c}
\xymatrix@C=20pt@R=20pt{ &\bullet\ar@/^0.25pc/[d]|{c}&\\
\bullet\ar@/^0.25pc/[r]|{b} &
\bullet\ar@/^0.25pc/[d]|{A}\ar@/^0.25pc/[u]|{C}\ar@/^0.25pc/[l]|{B}\ar@/^0.25pc/[r]|{D}&
\bullet\ar@/^0.25pc/[l]|{d} \\
&\star\ar@/^0.25pc/[u]|{a}&}\end{array}}} & {\small{
\begin{array}{c}
aA=0\quad cC=0\\bB=0\quad dD=0\\
Aa-Dd=Bb-Cc.
\end{array}}}
\end{array}
\]
We choose dimension vector and stability
\[
\alpha=\begin{array}{c}\xymatrix@C=-1pt@R=-1pt{&1&\\1&2&1\\ &1&}\end{array}\qquad \theta=\begin{array}{c}\xymatrix@C=-1pt@R=-1pt{&1&\\1&1&1\\ &\minus 5&}\end{array}
\]
and consider the resulting moduli space.  By choice of stability parameter, to specify an open set in this space we must
\begin{itemize}
\item specify, for each of the three vertices which are not either $\star$ or the middle vertex, a non-zero path (which we can change basis to assume to be the identity) from $\star$ to that vertex.
\item specify paths $(0\,\,1)$ and $(1\,\,0)$ from $\star$ to the middle vertex.
\end{itemize}
Different choices in the above lead to different open sets. Note
that we must be able to make such choices for any $\theta$-stable module
$M$ since by definition $M$ is $\star$-generated and so paths leaving
the trivial vertex must generate the vector spaces at all other
vertices.   For a stable $M$, it must be true that $a\neq 0$ and so after changing basis we can (and will) always assume  that $a=(1\,\,0)$.}

\t{Define the open sets $U_0$, $U_1$, $U_2$, $U_+$ and $U_-$ by the
following conditions: {\scriptsize{
\[
\begin{array}{ccccccc}
U_0 && aB=1& aC=1& aBbD=1&a=(1\,\,0)& b=(0\,\,1)\\
U_1 && aB=1& aC=1& aD=1&a=(1\,\,0)& b=(0\,\,1)\\
U_2 && aB=1& aC=1& aD=1&a=(1\,\,0)& d=(0\,\,1)\\
U_+ && aB=1& aDdC=1& aD=1&a=(1\,\,0)& d=(0\,\,1)\\
U_- && aDdB=1& aC=1& aD=1& a=(1\,\,0)& d=(0\,\,1)
\end{array}
\]
}}Pictorially we draw this as follows: {\scriptsize{
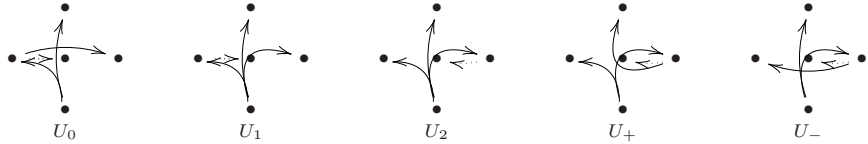
\[
\begin{array}{ccccc}
\begin{array}{c}
\xymatrix@C=10pt@R=10pt{ &\bullet&\\
\bullet\ar@{.>}[r]\ar@/^0.35pc/[0,2]& \bullet&
 \bullet \\
&\bullet\ar@<-0.25ex>@/^0.35pc/[-2,0]\ar@<0.75ex>@/_0.75pc/[-1,-1]&}\end{array}
& \begin{array}{c}
\xymatrix@C=10pt@R=10pt{ &\bullet&\\
\bullet\ar@{.>}[r]& \bullet&
 \bullet \\
&\bullet\ar@<-0.25ex>@/^0.35pc/[-2,0]\ar@<0.75ex>@/_0.75pc/[-1,-1]\ar@/^1.25pc/[-1,1]&}\end{array}
&\begin{array}{c}
\xymatrix@C=10pt@R=10pt{ &\bullet&\\
\bullet& \bullet&
 \bullet\ar@/^0.15pc/@{.>}[l] \\
&\bullet\ar@<-0.25ex>@/^0.35pc/[-2,0]\ar@<0.75ex>@/_0.75pc/[-1,-1]\ar@/^1.25pc/[-1,1]&}\end{array}
&\begin{array}{c}
\xymatrix@C=10pt@R=10pt{ &\bullet&\\
\bullet& \bullet&
 \bullet\ar@/^0.15pc/@{.>}[l]\ar@/^1.4pc/[-1,-1] \\
&\bullet\ar@<0.75ex>@/_0.75pc/[-1,-1]\ar@/^1.25pc/[-1,1]&}\end{array}
&\begin{array}{c}
\xymatrix@C=10pt@R=10pt{ &\bullet&\\
\bullet& \bullet&
 \bullet\ar@/^0.15pc/@{.>}[l]\ar@/^0.4pc/[0,-2] \\
&\bullet\ar@<-0.25ex>@/^0.35pc/[-2,0]\ar@/^1.25pc/[-1,1]&}\end{array}\\
U_0&U_1&U_2&U_+&U_-
\end{array}
\]
}}where the solid black lines correspond to the identity, and the
dotted arrow corresponds to the choice of vector $(0\,\,1)$.  To prove that these actually cover the moduli is perhaps the hardest part:}
\begin{lemma}
The open sets  $U_0$, $U_1$, $U_2$, $U_+$ and $U_-$ cover the moduli space.
\end{lemma}
\begin{proof}
Take an arbitrary stable module
{\scriptsize{
\[
M=\begin{array}{c} \xymatrix@C=30pt@R=30pt{
&{\C{}}\ar@/^0.5pc/[d]|(0.5){\qquad{\tiny{\left(\begin{smallmatrix}
c_1& c_2
\end{smallmatrix}\right)}} }&\\
{\C{}}\ar@/^0.5pc/[r]\ar@{}@<3ex>[r]|(0.4){{\tiny{\left(\begin{smallmatrix} b_1& b_2
\end{smallmatrix}\right) }}}&
{\C{2}}\ar@/^0.5pc/[d]|(0.5){\quad{\tiny{\left(\begin{smallmatrix}A_1\\
A_2\end{smallmatrix}\right)
}}}\ar@/^0.5pc/[u]|(0.5){{\tiny{\left(\begin{smallmatrix} C_1\\C_2
\end{smallmatrix}\right) }}\quad}\ar@/^0.5pc/[l]|(0.5){{\tiny{\left(\begin{smallmatrix}
B_1\\ B_2
\end{smallmatrix}\right) }}}\ar@/^0.5pc/[r]|(0.5){{\tiny{\left(\begin{smallmatrix} D_1\\
D_2 \end{smallmatrix}\right) }}}&
 {\C{}}\ar@/^0.5pc/[l]\ar@{}@<3ex>[l]|(0.4){{\tiny{\left(\begin{smallmatrix} d_1&d_2
\end{smallmatrix}\right) }}} \\
&{\C{}}\ar@/^0.5pc/[u]|(0.5){{\tiny{\left(\begin{smallmatrix} 1& 0
\end{smallmatrix}\right) }}\quad}&}
\end{array}
\]
}}We must show that we can change basis so that $M$ belongs to one of the above open sets.  We firstly prove the following claim: if $d$ is a linear multiple of $(1\,\,0)$ then we may choose $b=(0\,\,1)$.  To see this, since we must generate the module $M$ from $\star$, for every vertex $+$, $-$ and $2$ one of the paths corresponding to a generator of the corresponding CM module must be non-zero.  By assumption on $d$ we must be able to choose either $b=(0\,\,1)$ or $c=(0\,\,1)$ so that we generate the middle vertex.  If we can choose $b=(0\,\,1)$ we are done, hence suppose $b$ is a linear multiple of $(1\,\,0)$.  In this case $b=(\mu\,\,0)$ and $d=(\lambda\,\,0)$ for some $\lambda,\nu\in\C{}$ and so the preprojective relations force $\mu B_1=0$ and $\lambda D_1=0$.  Consequently the paths $aBb$ and $aDd$ are zero.  By inspection of the generators of the CM modules at vertices $+$, $-$ and $2$ this forces $aB\neq 0$, $aC\neq 0$ and $aD\neq 0$.  In particular $D_1\neq 0$ and $B_1\neq 0$, but this in turn forces $\lambda=\mu=0$ which is a contradiction since the preprojective relation at the middle vertex cannot hold.
 
 With this claim, the proof is now easy: Suppose first $d$ can be chosen as $(0\,\,1)$.  If $aB=0$ then by inspecting the generators of the CM module at $+$ we must have $aDdB\neq 0$, which forces $aCcB\neq 0$ by the preprojective relations.  Thus $M$ is in $U_-$.  Hence we may suppose that $aB\neq 0$.  Now if $aC=0$ then we must reach vertex $-$ by the other generator, so $aDdC\neq 0$ and so $M$ is in $U_+$.  Thus we may also assume that $aC\neq 0$.  Now if $aD=0$ then $D_1=0$ so $D_2\neq 0$ so that we generate.  But $D_2=dD=0$, a contradiction.  Hence necessarily $aD\neq 0$ and so $M$ is in $U_2$.
 
The above covers the case when we can choose $d=(0\,\,1)$ so suppose now that this is not the case, i.e. $d=(\lambda\,\,0)$ for some $\lambda\in\C{}$.  By the above claim we know we may choose $b=(0\,\,1)$.  If $aD\neq 0$ then $D_1\neq 0$ from which $dD=0$ gives us that $\lambda=0$ and so $d=0$.  Consequently to generate at vertex $+$ we need $aB\neq 0$ and to generate at vertex $-$ we need $aC\neq 0$.  Thus $M$ is in $U_1$.  Hence we may assume that $aD=0$, and so to generate at vertex $2$ requires that $aBbD\neq 0$ which forces $aCcD\neq 0$.  Consequently $M$ is in $U_0$.
\end{proof}
\t{Note that there are \emph{lots} of other open covers we could take.  We now compute the open sets above - we do the $U_0$ calculation in full and just summarize the others.  Any stable module in $U_0$ looks like
{\scriptsize
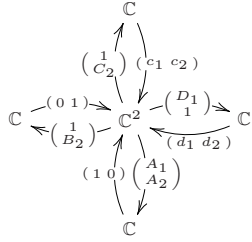
{
\[
\xymatrix@C=30pt@R=30pt{
&{\C{}}\ar@/^0.5pc/[d]|(0.5){\qquad{\tiny{\left(\begin{smallmatrix}
c_1& c_2
\end{smallmatrix}\right)}} }&&\\
{\C{}}\ar@/^0.5pc/[r]|(0.45){{\tiny{\left(\begin{smallmatrix} 0& 1
\end{smallmatrix}\right) }}}&
{\C{2}}\ar@/^0.5pc/[d]|(0.5){\quad{\tiny{\left(\begin{smallmatrix} A_1\\
A_2\end{smallmatrix}\right)
}}}\ar@/^0.5pc/[u]|(0.5){{\tiny{\left(\begin{smallmatrix} 1\\C_2
\end{smallmatrix}\right) }}\quad}\ar@/^0.5pc/[l]|(0.5){{\tiny{\left(\begin{smallmatrix}
1\\ B_2
\end{smallmatrix}\right) }}}\ar@/^0.5pc/[r]|(0.5){{\tiny{\left(\begin{smallmatrix} D_1\\
1 \end{smallmatrix}\right) }}}&
 {\C{}}\ar@/^0.5pc/[l]\ar@{}@<3ex>[l]|(0.4){{\tiny{\left(\begin{smallmatrix} d_1&d_2
\end{smallmatrix}\right) }}} \\
&{\C{}}\ar@/^0.5pc/[u]|(0.5){{\tiny{\left(\begin{smallmatrix} 1& 0
\end{smallmatrix}\right) }}\quad}&&}
\]
}}where the variables are scalars, subject only to the quiver
relations. Now {\small{
\[
\begin{array}{l}
   aA=0 \mbox{ implies } A_1=0\\
  bB=0 \mbox{ implies }B_2=0\\
  cC=0 \mbox{ implies }c_1=-c_2C_2\\
  dD=0 \mbox{ implies }d_2=-d_1D_1
\end{array}
\]
}}and so plugging this in our module becomes
{\scriptsize
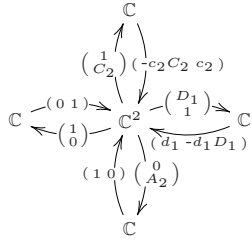
{
\[
\xymatrix@C=30pt@R=30pt{
&{\C{}}\ar@/^0.5pc/[d]|(0.5){\quad\qquad{\tiny{\left(\begin{smallmatrix}
\minus c_2C_2& c_2
\end{smallmatrix}\right)}} }&&\\
{\C{}}\ar@/^0.5pc/[r]|(0.45){{\tiny{\left(\begin{smallmatrix} 0& 1
\end{smallmatrix}\right) }}}&
{\C{2}}\ar@/^0.5pc/[d]|(0.5){\quad{\tiny{\left(\begin{smallmatrix} 0\\
A_2\end{smallmatrix}\right)
}}}\ar@/^0.5pc/[u]|(0.5){{\tiny{\left(\begin{smallmatrix} 1\\C_2
\end{smallmatrix}\right) }}\quad}\ar@/^0.5pc/[l]|(0.5){{\tiny{\left(\begin{smallmatrix}
1\\ 0
\end{smallmatrix}\right) }}}\ar@/^0.5pc/[r]|(0.5){{\tiny{\left(\begin{smallmatrix} D_1\\
1 \end{smallmatrix}\right) }}}&
 {\C{}}\ar@/^0.5pc/[l]\ar@{}@<3ex>[l]|(0.4){{\tiny{\left(\begin{smallmatrix} d_1&\minus d_1D_1
\end{smallmatrix}\right) }}}\\
&{\C{}}\ar@/^0.5pc/[u]|(0.5){{\tiny{\left(\begin{smallmatrix} 1& 0
\end{smallmatrix}\right) }}\quad}&&}
\]
}}But now there is only one relation left, namely $Aa-Dd=Bb-Cc$.
This gives {\scriptsize{
\[
\begin{pmatrix}
  0&0\\A_2&0
\end{pmatrix} -
\begin{pmatrix}
  d_1D_1&- d_1D_1^2\\d_1&- d_1D_1
\end{pmatrix}=
\begin{pmatrix}
  0&1\\0&0
\end{pmatrix} +
\begin{pmatrix}
  - c_2C_2&c_2\\  -c_2C_2^2&c_2C_2
\end{pmatrix}
\]
}}which yields the four conditions
 {\scriptsize{
\[
\begin{array}{c}
d_1D_1=c_2C_2\\
d_1D_1^2=1+ c_2\\
A_2=d_1-c_2C_2^2\\
d_1D_1=c_2C_2
\end{array}
\]
}}The second and third conditions eliminate the variables $c_2$ and
$A_2$, whereas the first and last conditions are the same.
Substituting the second condition into the first we see that this
open set is completely parameterized by $d_1$, $D_1$ and $C_2$
subject to the one relation $d_1D_1=(d_1D_1^2-1)C_2$, so $U_0$ is a smooth hypersurface in $\C{3}$.}
\t{Similarly we have {\scriptsize{
\[
\begin{array}{cccc}
U_1 & \begin{array}{c} \xymatrix@C=30pt@R=30pt{
&{\C{}}\ar@/^0.5pc/[d]|(0.5){\quad\qquad{\tiny{\left(\begin{smallmatrix}
\minus c_2C_2& c_2
\end{smallmatrix}\right)}} }&\\
{\C{}}\ar@/^0.5pc/[r]|(0.5){{\tiny{\left(\begin{smallmatrix} 0& 1
\end{smallmatrix}\right) }}}&
{\C{2}}\ar@/^0.5pc/[d]|(0.5){\quad{\tiny{\left(\begin{smallmatrix} 0\\
A_2\end{smallmatrix}\right)
}}}\ar@/^0.5pc/[u]|(0.5){{\tiny{\left(\begin{smallmatrix} 1\\C_2
\end{smallmatrix}\right) }}\quad}\ar@/^0.5pc/[l]|(0.5){{\tiny{\left(\begin{smallmatrix}
1\\ 0
\end{smallmatrix}\right) }}}\ar@/^0.5pc/[r]|(0.5){{\tiny{\left(\begin{smallmatrix} 1\\
D_2 \end{smallmatrix}\right) }}}&
 {\C{}}\ar@/^0.5pc/[l]\ar@{}@<3ex>[l]|(0.3){{\tiny{\left(\begin{smallmatrix} \minus d_2D_2&d_2
\end{smallmatrix}\right) }}} \\
&{\C{}}\ar@/^0.5pc/[u]|(0.5){{\tiny{\left(\begin{smallmatrix} 1& 0
\end{smallmatrix}\right) }}\quad}&}
\end{array}& \begin{array}{c}
c_2C_2=d_2D_2\\
 c_2=1+d_2\\
A_2=c_2C_2^2-d_2D_2^2\\
c_2C_2= d_1D_1
\end{array} & \C{3}_{d_2,D_2,C_2}/ (1+d_2)C_2=d_2D_2
\end{array}
\]
\[
\begin{array}{cccc}
U_2 & \begin{array}{c} \xymatrix@C=30pt@R=30pt{
&{\C{}}\ar@/^0.5pc/[d]|(0.5){\quad\qquad{\tiny{\left(\begin{smallmatrix}
\minus c_2C_2& c_2
\end{smallmatrix}\right)}} }&\\
{\C{}}\ar@/^0.5pc/[r]\ar@{}@<3ex>[r]|(0.4){{\tiny{\left(\begin{smallmatrix}
\minus b_2B_2&b_2
\end{smallmatrix}\right) }}}&
{\C{2}}\ar@/^0.5pc/[d]|(0.5){\quad{\tiny{\left(\begin{smallmatrix} 0\\
A_2\end{smallmatrix}\right)
}}}\ar@/^0.5pc/[u]|(0.5){{\tiny{\left(\begin{smallmatrix} 1\\C_2
\end{smallmatrix}\right) }}\quad}\ar@/^0.5pc/[l]|(0.5){{\tiny{\left(\begin{smallmatrix}
1\\ B_2
\end{smallmatrix}\right) }}}\ar@/^0.5pc/[r]|(0.5){{\tiny{\left(\begin{smallmatrix} 1\\
0 \end{smallmatrix}\right) }}}&
 {\C{}}\ar@/^0.5pc/[l]|(0.4){{\tiny{\left(\begin{smallmatrix} 0&1
\end{smallmatrix}\right) }}} \\
&{\C{}}\ar@/^0.5pc/[u]|(0.5){{\tiny{\left(\begin{smallmatrix} 1& 0
\end{smallmatrix}\right) }}\quad}&}
\end{array}& \begin{array}{c}
b_2B_2=c_2C_2\\
b_2=c_2-1\\
A_2=c_2C_2^2-b_2B_2^2\\
b_2B_2= c_2C_2
\end{array} & \C{3}_{b_2,B_2,C_2}/ (1+b_2)C_2=b_2B_2
\end{array}
\]
\[
\begin{array}{cccc}
U_+ & \begin{array}{c} \xymatrix@C=30pt@R=30pt{
&{\C{}}\ar@/^0.5pc/[d]|(0.5){\quad\qquad{\tiny{\left(\begin{smallmatrix}
c_1& \minus c_1C_1
\end{smallmatrix}\right)}} }&\\
{\C{}}\ar@/^0.5pc/[r]\ar@{}@<3ex>[r]|(0.4){{\tiny{\left(\begin{smallmatrix}
\minus b_2B_2&b_2
\end{smallmatrix}\right) }}}&
{\C{2}}\ar@/^0.5pc/[d]|(0.5){\quad{\tiny{\left(\begin{smallmatrix} 0\\
A_2\end{smallmatrix}\right)
}}}\ar@/^0.5pc/[u]|(0.5){{\tiny{\left(\begin{smallmatrix} C_1\\1
\end{smallmatrix}\right) }}\quad}\ar@/^0.5pc/[l]|(0.5){{\tiny{\left(\begin{smallmatrix}
1\\ B_2
\end{smallmatrix}\right) }}}\ar@/^0.5pc/[r]|(0.5){{\tiny{\left(\begin{smallmatrix} 1\\
0 \end{smallmatrix}\right) }}}&
 {\C{}}\ar@/^0.5pc/[l]|(0.4){{\tiny{\left(\begin{smallmatrix} 0&1
\end{smallmatrix}\right) }}} \\
&{\C{}}\ar@/^0.5pc/[u]|(0.5){{\tiny{\left(\begin{smallmatrix} 1& 0
\end{smallmatrix}\right) }}\quad}&}
\end{array}& \begin{array}{c}
b_2B_2=-c_1C_1\\
b_2=-c_1C_1^2-1\\
A_2=-b_2B_2^2-c_1\\
b_2B_2=-c_1C_1
\end{array} & \C{3}_{c_1,B_2,C_1}/ (1+c_1C_1^2)B_2=c_1C_1
\end{array}
\]
\[
\begin{array}{cccc}
U_- & \begin{array}{c} \xymatrix@C=30pt@R=30pt{
&{\C{}}\ar@/^0.5pc/[d]|(0.5){\quad\qquad{\tiny{\left(\begin{smallmatrix}
\minus c_2C_2& c_2
\end{smallmatrix}\right)}} }&\\
{\C{}}\ar@/^0.5pc/[r]\ar@{}@<3ex>[r]|(0.4){{\tiny{\left(\begin{smallmatrix}
b_1&\minus b_1B_1
\end{smallmatrix}\right) }}}&
{\C{2}}\ar@/^0.5pc/[d]|(0.5){\quad{\tiny{\left(\begin{smallmatrix} 0\\
A_2\end{smallmatrix}\right)
}}}\ar@/^0.5pc/[u]|(0.5){{\tiny{\left(\begin{smallmatrix} 1\\C_2
\end{smallmatrix}\right) }}\quad}\ar@/^0.5pc/[l]|(0.5){{\tiny{\left(\begin{smallmatrix}
B_1\\ 1
\end{smallmatrix}\right) }}}\ar@/^0.5pc/[r]|(0.5){{\tiny{\left(\begin{smallmatrix} 1\\
0 \end{smallmatrix}\right) }}}&
 {\C{}}\ar@/^0.5pc/[l]|(0.4){{\tiny{\left(\begin{smallmatrix} 0&1
\end{smallmatrix}\right) }}} \\
&{\C{}}\ar@/^0.5pc/[u]|(0.5){{\tiny{\left(\begin{smallmatrix} 1& 0
\end{smallmatrix}\right) }}\quad}&}
\end{array}& \begin{array}{c}
b_1B_1=-c_2C_2\\
c_2=1-b_1B_1^2\\
A_2=c_2C_2^2+b_1\\
b_1B_1=-c_2C_2
\end{array} & \C{3}_{b_1,B_1,C_2}/ (b_1B_1^2-1)C_2=b_1B_1
\end{array}
\]
}}Note in $U_2$ above the equation $b_2=c_2-1$ really means that we
have a choice of co-ordinate between $b_2$ and $c_2$; thus we could
equally well parameterize $U_2$ as $\C{3}_{c_2,B_2,C_2}/c_2C_2=(c_2-1)B_2$.}

\t{Hence we see that the space is covered by 5 open sets, each a smooth hypersurface in $\C{3}$.   By changing basis on the quiver it is also quite easy to write down the glues:%\footnote{The ordering of co-ordinates on RHS is as the ordering above.  Note we use the co-ordinates $(b_2,B_2,C_2)$ for $U_2$ on RHS.  Note also that in $U_1-U_2$ glue $\minus C_2=d_2(C_2-D_2)$ and so really "two" $d_2$'s!}
\[
\begin{array}{rcl} U_0\ni (d_1,D_1,C_2)&\longleftrightarrow&(\minus
d_1D_1^2,D_1^{-1},C_2)\in U_1\\ 
U_1\ni (d_2,D_2,C_2)&\longleftrightarrow&(d_2^{-1},\minus d_2D_2,\minus C_2)\in U_2\\ 
U_2\ni (c_2,B_2,C_2)&\longleftrightarrow&(\minus c_2C_2^2,B_2,C_2^{-1})\in U_+\\
U_2\ni (b_2,B_2,C_2)&\longleftrightarrow&(\minus b_2B_2^2,B_2^{-1},C_2)\in U_- 
\end{array}
\]
from which we can just see the configuration of $\mathbb{P}^1$'s.  The picture of the glues should (roughly) coincide with the picture drawn earlier.}
\end{example}
\begin{remark}
\t{Similar calculations to the above can be found in \cite{Leng} and \cite{Alvaro}.}
\end{remark}

We shall now illustrate how the reconstruction algebra changes this moduli picture:
\begin{example}
\t{Consider the group $\mathbb{D}_{5,3}$.  By Theorem~\ref{iso_of_rings1} the moduli presentation of the reconstruction algebra is
\[
\begin{array}{cc}
{\scriptsize{\begin{array}{c}
\xymatrix@C=20pt@R=20pt{ &\bullet\ar@/^0.25pc/[d]|{c}&\\
\bullet\ar@/^0.25pc/[r]|{b} &
\bullet\ar@/^0.25pc/[d]|{A}\ar@/^0.25pc/[u]|{C}\ar@/^0.25pc/[l]|{B}\ar@/^0.25pc/[r]|{D}&
\bullet\ar@/^0.25pc/[l]|{d}\ar@[green][1,-1]|{k_1} \\
&\star\ar@/^0.25pc/[u]|{a}&}\end{array}}} & {\scriptsize{
\begin{array}{c}
aA=0\quad cC=0\\bB=0\quad dD=0\\
Aa-Dd=Bb-Cc\\
k_1aD=dBbD\\
aDk_1=aCcA
\end{array}}}
\end{array}
\]
Choose dimension vector and stability as in the previous example.  Notice that \emph{the same} conditions that defined an open cover in the previous example give an open cover here, since arrows pointing to $\star$ don't add more choices.}

\t{Now the calculation from the previous example tells us almost everything, except now we have a new variable $k_1$ inside every open set.  The point is that the \emph{only} open set which changes is $U_0$.  The reason for this is quite simple: in the relations of the reconstruction algebra $k_1aD=dBbD$ and notice that $aD=1$ in every open set except $U_0$.  Thus $k_1=dBbD$ in every open set except $U_0$ and consequently we can put $k_1$ in terms of the other variables.  Hence $k_1$ isn't really an extra variable in these open sets, so they do not change.}   

\t{Thus the only open set that changes is $U_0$, and by the previous calculation (the relations above the line) and our new relations (shown below the line) we see that $U_0$ is given by
{\scriptsize{
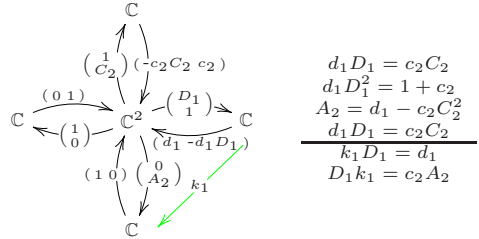
\[
\begin{array}{cc}
\begin{array}{c} \xymatrix@C=30pt@R=30pt{
&{\C{}}\ar@/^0.5pc/[d]|(0.5){\quad\qquad{\tiny{\left(\begin{smallmatrix}
\minus c_2C_2& c_2
\end{smallmatrix}\right)}} }&\\
{\C{}}\ar@/^0.5pc/[r]\ar@{}@<3ex>[r]|(0.4){{\tiny{\left(\begin{smallmatrix}
0&1
\end{smallmatrix}\right) }}}&
{\C{2}}\ar@/^0.5pc/[d]|(0.5){\quad{\tiny{\left(\begin{smallmatrix} 0\\
A_2\end{smallmatrix}\right)
}}}\ar@/^0.5pc/[u]|(0.5){{\tiny{\left(\begin{smallmatrix} 1\\C_2
\end{smallmatrix}\right) }}\quad}\ar@/^0.5pc/[l]|(0.5){{\tiny{\left(\begin{smallmatrix}
1\\ 0
\end{smallmatrix}\right) }}}\ar@/^0.5pc/[r]|(0.5){{\tiny{\left(\begin{smallmatrix} D_1\\
1 \end{smallmatrix}\right) }}}&
 {\C{}}\ar@/^0.5pc/[l]\ar@{}@<3ex>[l]|(0.4){{\tiny{\left(\begin{smallmatrix} d_1&\minus d_1D_1
\end{smallmatrix}\right) }}}\ar@[green]@<2ex>[1,-1]|{k_1} \\
&{\C{}}\ar@/^0.5pc/[u]|(0.5){{\tiny{\left(\begin{smallmatrix} 1& 0
\end{smallmatrix}\right) }}\quad}&}
\end{array}& \begin{array}{c}
d_1D_1=c_2C_2\\
d_1D_1^2=1+c_2\\
A_2=d_1-c_2C_2^2\\
d_1D_1=c_2C_2\\ \hline
k_1D_1=d_1\\
D_1k_1=c_2A_2
\end{array}
\end{array}
\]
}}Since $d_1=k_1D_1$, instead of being given by $d_1,D_1,C_2$ subject to $d_1D_1=(d_1D_1^2-1)C_2$ the open set $U_0$ is now given by $k_1,D_1,C_2$ subject to $k_1D_1^2=(k_1D_1^3-1)C_2$  .  Also, the gluing between $U_0$ and $U_1$ has changed to
\[
\xymatrix@C=30pt{U_0\ni(k_1,D_1,C_2)\ar@{<->}[r]^(0.3){D_1\neq 0}&(\minus
(k_1D_1)D_1^2,D_1^{-1},C_2)=(\minus k_1D_1^3,D_1^{-1},C_2)\in U_1 }.
\]
Thus we see that the red curve has changed into a $(-3)$-curve, nothing else in the open cover has changed and so the dual graph is now} 
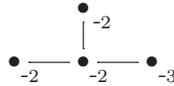
\[
\xymatrix@C=15pt@R=10pt{ &\bullet\ar@{-}[d]^<{\minus 2}&&\\
\bullet\ar@{-}[r]_<{\minus 2} & \bullet\ar@{-}[r]_<{\minus 2}
& \bullet\ar@{}[r]_<{\minus 3}&}
\]
\end{example}


\begin{thebibliography}{Wem09}
\bibitem[BR78]{Riemenschneider_dihedral}
K.~Behnke and O.~Riemenschneider, \emph{Diedersingularit{\"{a}}ten},
Abh. Math. Sem. Univ. Hamburg \textbf{47} (1978), 210--227.

\bibitem[BSW08]{BSW}
R.~Bocklandt, T.~Schedler and M.~Wemyss, \emph{Superpotentials and Higher Order Derivations}, arXiv:0802.0162 (2008)

\bibitem[Bri02]{Bridgeland}
T.~Bridgeland, \emph{Flops and derived categories}, Invent. Math. \textbf{147} (2002), no.~3, 613--632.

\bibitem[Bri68]{Brieskorn}
E.~Brieskorn, \emph{Rationale singularit{\"{a}}ten komplexer fl{\"{a}}chen},
  Invent. Math. \textbf{4} (1968), 336--358.

\bibitem[BIRS]{BIRS}
A. Buan, O. Iyama, I. Reiten and D. Smith, \emph{Mutation of cluster-tilting objects and potentials},
arXiv:0804.3813 (2008)

\bibitem[CB98]{CBH}
W.~Crawley-Boevey and M.~P. Holland, \emph{Noncommutative deformations
  of {K}leinian singularities}, Duke Math. J. \textbf{92} (1998), no.~3,
  605--635.

\bibitem[Ish02]{Ishii} 
A.~Ishii, \emph{On the McKay correspondence for a finite small subgroup of ${\rm GL}(2,\Bbb C)$}, J. Reine Angew. Math. \textbf{549} (2002), 221--233

\bibitem[IW08]{Iyama_Wemyss_specials}
O.~Iyama and M.~Wemyss, \emph{The classification of special Cohen Macaulay modules},
arXiv:0809.1958 (2008)

\bibitem[Len02]{Leng}
R.~Leng, \emph{The McKay Correspondence and Orbifold Riemann-Roch}, Warwick PhD thesis (2002).

\bibitem[NdC09]{Alvaro}
A.~Nolla de Celis, \emph{Dihedral $G$-Hilb via representations of the McKay quiver}, arXiv:0905.1195 (2009).

\bibitem[Rie77]{Riemenschneider_invarianten}
O.~Riemenschneider, \emph{Invarianten endlicher {U}ntergruppen},
Math. Z
  \textbf{153} (1977), 37--50.

\bibitem[Wem07]{Wemyss_reconstruct_A}
M.~Wemyss, \emph{Reconstruction algebras of type ${A}$},
arXiv:0704.3693 (2007).
  
\bibitem[Wem08]{Wemyss_GL2}
\bysame, \emph{The ${GL}(2)$ {M}c{K}ay correspondence},
arXiv:0809.1973 (2008).

\bibitem[Wem09]{Wemyss_reconstruct_D(i)}
\bysame, \emph{Reconstruction algebras of type ${D}$ ({I})}, arXiv:0905.1154 (2009).

\bibitem[Wun88]{Wunram_generalpaper}
J.~Wunram, \emph{Reflexive modules on quotient surface singularities},
  Mathematische Annalen \textbf{279} (1988), no.~4, 583--598.


\end{thebibliography}
\end{document}